\documentclass[reqno]{amsart}

\newtheorem{Theorem}{Theorem}[section]
\newtheorem{Proposition}[Theorem]{Proposition}
\newtheorem{Corollary}[Theorem]{Corollary}
\newtheorem{Lemma}[Theorem]{Lemma}
\theoremstyle{remark}
\newtheorem{Remark}[Theorem]{Remark}
\numberwithin{equation}{section}

\allowdisplaybreaks

\def\bmcr{\hbox{$\boldsymbol C_{\boldsymbol r}$}}
\def\bmar{\hbox{$\boldsymbol A_{\boldsymbol {r-1}}$}}
\def\bsixphifive{\hbox{$ _{\boldsymbol 6}\boldsymbol\phi_{\boldsymbol 5}$}}
\def\bthreephitwo{\hbox{$ _{\boldsymbol 3}\boldsymbol\phi_{\boldsymbol 2}$}}
\def\btwophione{\hbox{$ _{\boldsymbol 2}\boldsymbol\phi_{\boldsymbol 1}$}}
\def\bonephizero{\hbox{$ _{\boldsymbol 1}\boldsymbol\phi_{\boldsymbol 0}$}}

\begin{document}

\title[Summations and transformations]
{Summations and transformations for\\
multiple basic and elliptic hypergeometric\\
 series by determinant evaluations}

\author{Hjalmar Rosengren}
\address{Department of Mathematics, Chalmers University of Technology and
G\"oteborg University, SE-412 96 G\"oteborg, Sweden}
\email{hjalmar@math.chalmers.se}
\urladdr{http://www.math.chalmers.se/{\textasciitilde}hjalmar}

\author[Michael Schlosser]{Michael Schlosser$^*$}
\address{Institut f\"ur Mathematik der Universit\"at Wien,
Strudlhofgasse 4, A-1090 Wien, Austria}
\email{schlosse@ap.univie.ac.at}
\urladdr{http://www.mat.univie.ac.at/{\textasciitilde}schlosse}

\thanks{$^*$The second author was supported by an APART grant of the Austrian
Academy of Sciences}
\date{first version: April 17, 2003; revised: June 16, 2003}

\dedicatory{Dedicated to Tom Koornwinder}

\begin{abstract}
Using multiple $q$-integrals and a determinant evaluation,
we establish a multivariable extension of Bailey's
nonterminating $_{10}\phi_9$ transformation.
From this result, we deduce new multivariable
terminating $_{10}\phi_9$ transformations, 
$_8\phi_7$ summations and other identities. We also use similar methods to
derive new multivariable ${}_1\psi_1$ summations.
Some of our results are extended to the case of elliptic
hypergeometric series.
\end{abstract}

\maketitle

\section{Introduction}\label{sec0}

Basic hypergeometric series 
have various applications in combinatorics, number theory,
representation theory, statistics, and physics,
see Andrews~\cite{andappl}, \cite{qandrews}.  For a general account
of the importance of basic hypergeometric series in the theory of
special functions see Andrews, Askey, and Roy~\cite{sfaar}.

There are different types of multivariable basic hypergeometric
series extending the classical one-dimensional theory~\cite{grhyp}.
Several recent multivariable extensions can be associated to
{\em root systems} or, equivalently, to {\em Lie algebras}.
The specific series we consider in this paper have the advantage
that they can be conveniently studied from a purely analytical
point of view. In this respect, we understand the root system
terminology used in this paper (such as ``$A_{r-1}$ series'', or
``$C_r$ series'') simply as a classification of
certain multiple series according
to specific factors (such as a Vandermonde determinant)
appearing in the summand. We omit giving a precise definition here,
but refer instead to papers of Bhatnagar~\cite{bhatdn} or
Milne~\cite[Sec.~5]{milnetf}. We mention that, although these series
first arose (in the limit $q\rightarrow 1$) in the representation theory
of compact Lie groups \cite{hbl}, many questions remain about 
this connection. In particular, there is no known (quantum) group
interpretation of the type of series that we will study.

In this paper, we give a multivariable nonterminating $_{10}\phi_9$
transformation for the root system $C_r$ (or, equivalently,
the symplectic group $Sp(r)$), see Corollary~\ref{cnnt109}.
Our result extends a $C_r$ nonterminating $_8\phi_7$ summation
recently found by one of us in \cite{schlcnnt87}. To our knowledge,
our new transformation formula is the first multivariable generalization
of Bailey's nonterminating $_{10}\phi_9$ transformation (see
\eqref{109ntgl} below) that has appeared in the literature.
We deduce this result from an equivalent multiple $q$-integral
transformation, Theorem~\ref{cnint109}. 
In our proof of the latter we utilize a simple determinant method,
essentially the same which was introduced by Gustafson and
Krattenthaler~\cite{guskradet} and which was further exploited
by one of us in \cite{schlhypdet} and \cite{schlcnnt87} to derive
a number of identities for multiple basic hypergeometric series.
This method was also recently
employed by Spiridonov [22, Th.~3] who used it to
derive a $C_r$ elliptic Nasrallah--Rahman integral.
Here, as in \cite{schlcnnt87}, we apply the determinant method to
{\em $q$-integrals} to derive our main result.
We also derive some new multivariable extensions of
Ramanujan's ${}_1\psi_1$ summation.

In the final section we briefly
 discuss \emph{elliptic} extensions of some of our identities. 
Elliptic (or modular) hypergeometric series is a recently introduced
 extension of basic hypergeometric series, which was
motivated by certain models in statistical mechanics \cite{ft}. 
 Warnaar \cite{warnell} used the determinant method to give
 an elliptic analogue of a $C_r$ Jackson summation 
  from  \cite{schlhypdet}. This identity was used by one of us \cite{r1}
to prove a second elliptic $C_r$ Jackson summation conjectured by Warnaar
 (the basic case of this identity was proved by van Diejen
 and Spiridonov \cite{ds1}; cf.\ \cite{ra} for a third proof), as well
 as  to prove an elliptic analogue of a third $C_r$
 Jackson  summation due to
 Denis and Gustafson \cite{dg} and Milne and Newcomb \cite{mn}; cf.\ \cite{r2}.
(Spiridonov and van Diejen showed that the second and third 
 summation also  follow from certain conjectured multiple integral
evaluations \cite{ds1}, \cite{ds2}.) Note that although all three Jackson
summations are connected to the root system $C_r$, they are 
different in nature. We want to stress here that
the first summation, which is related to the present work, 
is apparently the simplest but may be used to prove the other two.  
We hope that our new identities
will also be useful to study different, apparently more
complicated, types of multivariable series.

Our paper is organized as follows: In Section~\ref{secpre},
we review some basics in the theory of basic hypergeometric series.
We also note a determinant lemma which we need as an
ingredient in proving our results.
In Section~\ref{secmint}, we derive a multiple $q$-integral transformation,
Theorem~\ref{cnint109},
which in Section~\ref{secm109} is used to explicitly write out a
nonterminating $_{10}\phi_9$ transformation for the root system $C_r$,
see Corollary~\ref{cnnt109}. In Section 5 we specialize this identity
to obtain new terminating $C_r$ $_{10}\phi_9$ transformations. 
In Section 6 we derive  new multivariable extensions of
Ramanujan's ${}_1\psi_1$ summation. In Section 7 we give a
number of special and limit cases of our multivariable
$_{10}\phi_9$ transformations. Finally,
in Section 8 we prove that our terminating $C_r$ $_{10}\phi_9$
transformations and $_8\phi_7$ summations extend to the case of
elliptic hypergeometric series.

\section{Basic hypergeometric series and a determinant lemma}\label{secpre}

Here we recall some standard notation for
 basic hypergeometric series (cf.\ \cite{grhyp}).

Let $q$ be a complex number (called the ``base") such that $0<|q|<1$.
We define the
{\em $q$-shifted factorial} for all integers $k$ by
\begin{equation*}
(a;q)_{\infty}:=\prod_{j=0}^{\infty}(1-aq^j)\qquad
\text{and}\qquad
(a;q)_k:=\frac{(a;q)_{\infty}}{(aq^k;q)_{\infty}}.
\end{equation*}
Since we are working with the same base $q$ throughout this article,
we omit writing it out explicitly, i.e., we use
\begin{equation}\label{defqsf}
(a)_k:=(a;q)_k,
\end{equation}
where $k$ is an integer or infinity.
Further, we  employ the condensed notation
\begin{equation*}
(a_1,\ldots,a_m)_k\equiv (a_1)_k\dots(a_m)_k,
\end{equation*}
where $k$ is an integer or infinity.
We denote the {\em basic hypergeometric ${}_s\phi_{s-1}$ series} by
\begin{equation}\label{defhyp}
{}_s\phi_{s-1}\!\left[\begin{matrix}a_1,a_2,\dots,a_s\\
b_1,b_2,\dots,b_{s-1}\end{matrix}\,;q,z\right]:=
\sum _{k=0} ^{\infty}\frac {(a_1,a_2,\dots,a_s)_k}
{(q,b_1,\dots,b_{s-1})_k}z^k,
\end{equation}
and the {\em bilateral basic hypergeometric ${}_s\psi_s$ series} by
\begin{equation}\label{defhypb}
{}_s\psi_s\!\left[\begin{matrix}a_1,a_2,\dots,a_s\\
b_1,b_2,\dots,b_s\end{matrix}\,;q,z\right]:=
\sum _{k=-\infty} ^{\infty}\frac {(a_1,a_2,\dots,a_s)_k}
{(b_1,b_2,\dots,b_s)_k}z^k,
\end{equation}
respectively.  
See \cite[p.~25 and p.~125]{grhyp} for the criteria
of when these series terminate, or, if not, when they converge. 

The classical theory of basic hypergeometric series contains
numerous summation and transformation formulae
involving ${}_s\phi_{s-1}$ or ${}_s\psi_s$  series.
Many of these require that the parameters satisfy the condition of being
either balanced and/or very-well-poised.
An ${}_s\phi_{s-1}$ basic hypergeometric series is called
{\em balanced} if $b_1\cdots b_{s-1}=a_1\cdots a_sq$ and $z=q$.
An ${}_s\phi_{s-1}$ series is {\em well-poised} if
$a_1q=a_2b_1=\cdots=a_sb_{s-1}$ and {\em very-well-poised}
if it is well-poised and if $a_2=-a_3=q\sqrt{a_1}$.
Note that the factor
\begin{equation*}
\frac {1-a_1q^{2k}}{1-a_1}
\end{equation*}
appears in  a very-well-poised series.
Similarly, a bilateral ${}_s\psi_s$ basic hypergeometric series is
well-poised if $a_1b_1=a_2b_2\cdots=a_sb_s$ and very-well-poised if,
in addition, $a_1=-a_2=qb_1=-qb_2$.

The standard reference for basic hypergeometric series
is Gasper and Rahman's text~\cite{grhyp}.
In our computations throughout this paper we
frequently use some elementary identities of
$q$-shifted factorials, listed in \cite[Appendix~I]{grhyp}.

In the following we display some identities which we utilize
in the subsequent sections.

Bailey's~\cite[Eq.~(7.2)]{bailwp} nonterminating
very-well-poised $_{10}\phi_9$ summation,
\begin{multline}\label{109ntgl}
{}_{10}\phi_9\!\left[\begin{matrix}a,\,qa^{\frac 12},-qa^{\frac 12},
b,c,d,e,f,g,h\\a^{\frac 12},-a^{\frac 12},
aq/b,aq/c,aq/d,aq/e,aq/f,aq/g,aq/h\end{matrix}\,;q,q\right]\\
+\frac{(aq,b/a,c,d,e,f,g,h)_{\infty}}
{(b^2q/a,a/b,aq/c,aq/d,aq/e,aq/f,aq/g,aq/h)_{\infty}}\\\times
\frac{(bq/c,bq/d,bq/e,bq/f,bq/g,bq/h)_{\infty}}
{(bc/a,bd/a,be/a,bf/a,bg/a,bh/a)_{\infty}}\\\times
{}_{10}\phi_9\!\left[\begin{matrix}b^2/a,\,qba^{-\frac 12},-qba^{-\frac 12},
b,bc/a,bd/a,be/a,bf/a,bg/a,bh/a\\
ba^{-\frac 12},-ba^{-\frac 12},
bq/a,bq/c,bq/d,bq/e,bq/f,bq/g,bq/h\end{matrix}\,;q,q\right]\\
=\frac{(aq,b/a,\lambda q/f,\lambda q/g,\lambda q/h,
bf/\lambda,bg/\lambda,bh/\lambda)_{\infty}}
{(\lambda q,b/\lambda,aq/f,aq/g,aq/h,bf/a,bg/a,bh/a)_{\infty}}\\\times
{}_{10}\phi_9\!\left[\begin{matrix}\lambda,\,q\lambda^{\frac 12},
-q\lambda^{\frac 12},b,\lambda c/a,\lambda d/a,\lambda e/a,f,g,h\\
\lambda^{\frac 12},-\lambda^{\frac 12},\lambda q/b,
aq/c,aq/d,aq/e,\lambda q/f,\lambda q/g,\lambda q/h\end{matrix}\,;q,q\right]\\
+\frac{(aq,b/a,f,g,h,bq/f,bq/g,bq/h)_{\infty}}
{(b^2q/\lambda,\lambda/b,aq/c,aq/d,aq/e,aq/f,aq/g,aq/h)_{\infty}}\\\times
\frac{(\lambda c/a,\lambda d/a,\lambda e/a,
abq/\lambda c,abq/\lambda d,abq/\lambda e)_{\infty}}
{(bc/a,bd/a,be/a,bf/a,bg/a,bh/a)_{\infty}}\\\times
{}_{10}\phi_9\!\left[\begin{matrix}b^2/\lambda,\,qb\lambda^{-\frac 12},
-qb\lambda^{-\frac 12},b,bc/a,bd/a,be/a,bf/\lambda,bg/\lambda,bh/\lambda\\
b\lambda^{-\frac 12},-b\lambda^{-\frac 12},bq/\lambda,
abq/c\lambda,abq/d\lambda,abq/e\lambda,bq/f,bq/g,bq/h\end{matrix}\,;q,q\right],
\end{multline}
where $\lambda=a^2q/cde$ and $a^3q^2=bcdefgh$ (cf.~\cite[Eq.~(2.12.9)]{grhyp}),
stands on the top of the classical hierarchy of transformation formulae for
basic hypergeometric series.
If $h=q^{-n}$, where $n=0,1,2,\dots$, it
reduces to Bailey's terminating very-well-poised $_{10}\phi_9$
transformation
(cf.\ \cite[Eq.~(2.9.1)]{grhyp}):
\begin{multline}\label{109gl}
{}_{10}\phi_9\!\left[\begin{matrix}a,\,qa^{\frac 12},-qa^{\frac 12},
b,c,d,e,f,g,q^{-n}\\a^{\frac 12},-a^{\frac 12},
aq/b,aq/c,aq/d,aq/e,aq/f,aq/g,aq^{1+n}\end{matrix}\,;q,q\right]\\
=\frac {(aq,aq/ef,\lambda q/e,\lambda q/f)_n}
{(aq/e,aq/f,\lambda q/ef,\lambda q)_n}\\\times
{}_{10}\phi_9\!\left[\begin{matrix}\lambda,\,q\lambda^{\frac 12},
-q\lambda^{\frac 12},\lambda b/a,\lambda c/a,\lambda d/a,e,f,g,q^{-n}\\
\lambda^{\frac 12},-\lambda^{\frac 12},aq/b,aq/c,aq/d,
\lambda q/e,\lambda q/f,\lambda q/g,\lambda q^{1+n}\end{matrix}\,;q,q\right],
\end{multline}
where $\lambda=a^2q/bcd$ and $a^3q^{3+n}=bcdefg$.
For $cd= aq$ (hence $b\lambda= a$) Bailey's transformation
reduces to Jackson's~\cite{jacksum} terminating very-well-poised balanced
${}_8\phi_7$ summation which, after substitution of variables
($f\mapsto c$, $g\mapsto d$) can be written as
(cf.\ \cite[Eq.~(2.6.2)]{grhyp})
\begin{multline}\label{87gl}
{}_8\phi_7\!\left[\begin{matrix}a,\,qa^{\frac 12},-qa^{\frac 12},
b,c,d,e,q^{-n}\\a^{\frac 12},-a^{\frac 12},
aq/b,aq/c,aq/d,aq/e,aq^{1+n}\end{matrix}\,;q,q\right]\\
=\frac {(aq,aq/bc,aq/bd,aq/cd)_n}
{(aq/b,aq/c,aq/d,aq/bcd)_n},
\end{multline}
where $a^2q^{1+n}=bcde$.
 On the other hand, if $cd=aq$ (i.e., $e\lambda=a$)
\eqref{109ntgl} reduces to Bailey's nonterminating two-term
$_8\phi_7$ summation. These identities contain many other important
summation and transformation formulae for basic hypergeometric series.
For a sequence of derivations leading up to the nonterminating
$_{10}\phi_9$ transformation, see the exposition in
Gasper and Rahman~\cite[Section 2]{grhyp}.

For studying nonterminating basic hypergeometric series it is often
convenient to utilize Jackson's~\cite{jackson} $q$-integral notation,
defined by
\begin{equation}\label{qint}
\int_a^bf(t)\,d_qt=\int_0^bf(t)\;d_qt-\int_0^af(t)\;d_qt,
\end{equation}
where
\begin{equation}\label{qint0a}
\int_0^af(t)\;d_qt=a(1-q)\sum_{k=0}^{\infty}f(aq^k)q^k.
\end{equation}
If $f$ is continuous on $[0,a]$, then it is easily seen that
\begin{equation*}
\lim_{q\to1^-}\int_0^a f(t)\;d_qt=\int_0^a f(t)\;dt,
\end{equation*}
see~\cite[Eq.~(1.11.6)]{grhyp}.

Using the above $q$-integral notation, the nonterminating $_{10}\phi_9$
transformation \eqref{109ntgl} can be expressed as
\begin{multline}\label{int109gl}
\int_a^b\frac{(1-t^2/a)(qt/a,qt/b,qt/c,qt/d,qt/e,qt/f,qt/g,qt/h)_{\infty}}
{(t,bt/a,ct/a,dt/a,et/a,ft/a,gt/a,ht/a)_{\infty}}\;d_qt\\
=\frac a{\lambda}
\frac{(b/a,aq/b,\lambda c/a,\lambda d/a,\lambda e/a,
bf/\lambda,bg/\lambda,bh/\lambda)_{\infty}}
{(b/\lambda,\lambda q/b,c,d,e,bf/a,bg/a,bh/a)_{\infty}}\\\times
\int_{\lambda}^b\frac{(1-t^2/\lambda)
(qt/\lambda,qt/b,aqt/c\lambda,aqt/d\lambda,aqt/e\lambda,
qt/f,qt/g,qt/h)_{\infty}}
{(t,bt/\lambda,ct/a,dt/a,et/a,
ft/\lambda,gt/\lambda,ht/\lambda)_{\infty}}\;d_qt,
\end{multline}
where $\lambda=a^2q/cde$ and $a^3q^2=bcdefgh$
(cf.\ \cite[Eq.~(2.12.10)]{grhyp}).

One of the most important summmation theorems for 
bilateral basic hypergeometric series is
Ramanujan's ${}_1\psi_1$ summation~\cite[Eq.~(12.12.2)]{hardy},
which reads
\begin{equation}\label{11gl}
{}_1\psi_1\!\left[\begin{matrix}a\\b\end{matrix}\,;q,z\right]
=\frac {(q,az,q/az,b/a)_{\infty}}{(b,z,b/az,q/a)_{\infty}},
\end{equation}
where $|b/a|<|z|<1$ (cf.\ \cite[Eq.~(5.2.1)]{grhyp}).
This beautiful formula contains several important
identities as special cases, see Gasper and Rahman~\cite{grhyp}.

We conclude this section with a determinant evaluation, which will be
our main tool. It was given explicitly in \cite[Lemma~A.1]{schlhypdet},
as a special case of a determinant lemma
by Krattenthaler~\cite[Lemma~34]{krattmaj}.

\begin{Lemma}\label{lemdet1}
Let $X_1,\dots,X_r$, $A$, $B$, and $C$ be indeterminate. Then there holds
\begin{multline}\label{lemdet1gl}
\det_{1\le i,j\le r}\left(
\frac{(AX_i,AC/X_i)_{r-j}}
{(BX_i,BC/X_i)_{r-j}}\right)=
\prod_{1\le i<j\le r}(X_j-X_i)(1-C/X_iX_j)\\\times
A^{\binom r2}q^{\binom r3}
\prod_{i=1}^{r}\frac {(B/A,ABCq^{2r-2i})_{i-1}}
{(BX_i,BC/X_i)_{r-1}}.
\end{multline}
\end{Lemma}

The above determinant evaluation was generalized to the elliptic case
(more precisely, to an evaluation involving Jacobi theta functions)
by Warnaar~\cite[Cor.~5.4]{warnell}. We make use of the elliptic
version of Lemma~\ref{lemdet1} in Section~\ref{secell}.

\section{A multiple $q$-integral transformation}\label{secmint}

By iteration, the extension of \eqref{qint0a} to {\em multiple}
$q$-integrals is straightforward:
\begin{multline}\label{mqint0a}
\int_0^{a_1}\dots\int_0^{a_r}f(t_1,\dots,t_r)\;d_qt_r\dots d_qt_1\\
=a_1\dots a_r(1-q)^r\sum_{k_1,\dots,k_r=0}^{\infty}f(a_1q^{k_1},\dots,
a_rq^{k_r})q^{k_1+\dots+k_r}.
\end{multline}
Similarly, the extension of \eqref{qint} is
\begin{multline}\label{mqint}
\int_{a_1}^{b_1}\dots\int_{a_r}^{b_r}f(t_1,\dots,t_r)\;d_qt_r\dots d_qt_1\\
=\sum_{S\subseteq\{1,2,\dots,r\}}\left(\prod_{i\in S}(-a_i)\right)
\left(\prod_{i\notin S}b_i\right)(1-q)^r\\\times
\sum_{k_1,\dots,k_r=0}^{\infty}f(c_1(S)q^{k_1},\dots,c_r(S)q^{k_r})
q^{k_1+\dots+k_r},
\end{multline}
where the outer sum runs over all $2^r$ subsets $S$ of $\{1,2,\dots,r\}$,
and where $c_i(S)=a_i$ if $i\in S$ and $c_i(S)=b_i$ if $i\notin S$,
for $i=1,\dots,r$.

We give our main result, a $C_r$ extension of \eqref{int109gl}:

\begin{Theorem}\label{cnint109}
Let $a^3q^{3-r}=bc_id_ie_ix_ifgh$ and $\lambda=a^2q/c_id_ie_ix_i$
 for $i=1,\dots,r$.
Then there holds
\begin{multline}\label{cnint109gl}
\int_{ax_1}^{b}\dots\int_{ax_r}^{b}
\prod_{1\le i<j\le r}(t_i-t_j)(1-t_it_j/a)
\prod_{i=1}^r(1-t_i^2/a)\frac{(qt_i/ax_i,qt_i/b)_{\infty}}
{(t_ix_i,bt_i/a)_{\infty}}\\
\times\prod_{i=1}^r
\frac{(qt_i/c_i,qt_i/d_i,qt_i/e_i,qt_i/f,qt_i/g,qt_i/h)_{\infty}}
{(c_it_i/a,d_it_i/a,e_it_i/a,ft_i/a,gt_i/a,ht_i/a)_{\infty}}\;
d_qt_r\dots d_qt_1\\
=\left(\frac a\lambda\right)^{\binom{r+1}2}
\prod_{i=1}^r\frac
{(b/ax_i,ax_iq/b,\lambda c_ix_i/a,\lambda d_ix_i/a,\lambda e_ix_i/a)_\infty}
{(b/\lambda x_i,\lambda x_iq/b,c_ix_i,d_ix_i,e_ix_i)_\infty}\\\times
\prod_{i=1}^r\frac
{(bfq^{i-1}/\lambda,bgq^{i-1}/\lambda,bhq^{i-1}/\lambda)_\infty}
{(bfq^{i-1}/a,bgq^{i-1}/a, bhq^{i-1}/a)_\infty}\\\times
\int_{\lambda x_1}^{b}\dots\int_{\lambda x_r}^{b}
\prod_{1\le i<j\le r}(t_i-t_j)(1-t_it_j/\lambda)
\prod_{i=1}^r(1-t_i^2/\lambda)\frac{(qt_i/\lambda x_i,qt_i/b)_{\infty}}
{(t_ix_i,bt_i/\lambda)_{\infty}}\\
\times\prod_{i=1}^r
\frac{(aqt_i/c_i\lambda,aqt_i/d_i\lambda,aqt_i/e_i\lambda,qt_i/f,
qt_i/g,qt_i/h)_{\infty}}
{(c_it_i/a,d_it_i/a,e_it_i/a,ft_i/\lambda,
gt_i/\lambda,ht_i/\lambda)_{\infty}}\;d_qt_r\dots d_qt_1.
\end{multline}
\end{Theorem}

\begin{proof}
We have
\begin{multline*}
\prod_{1\le i<j\le r}(t_i-t_j)(1-t_it_j/a)
=\prod_{i=1}^r\frac{(q^{2-r}t_i/g,aq^{2-r}/gt_i)_{r-1}}
{(aq^{2-r}/fg,fq^{2+r-2i}/g)_{i-1}}t_i^{r-1}\\\times
f^{-\binom r2}q^{-\binom r3}\;
\det_{1\le i,j\le r}\left(\frac {(ft_i/a,f/t_i)_{r-j}}
{(q^{2-r}t_i/g,aq^{2-r}/gt_i)_{r-j}}\right),
\end{multline*}
due to the $X_i\mapsto t_i$, $A\mapsto f/a$, $B\mapsto q^{2-r}/g$,
and $C\mapsto a$ case of Lemma~\ref{lemdet1}.
Hence, using some elementary identities from \cite[Appendix~I]{grhyp},
we may write the left-hand side of \eqref{cnint109gl} as
\begin{multline*}
\left(\frac ag\right)^{\binom r2}q^{-\binom r3}
\prod_{i=1}^r(aq^{2-r}/fg,fq^{2+r-2i}/g)_{i-1}^{-1}\\\times
\det_{1\le i,j\le r}\Bigg(\int_{ax_i}^b
\frac {(1-t_i^2/a)(t_iq/ax_i,t_iq/b)_{\infty}}
{(t_ix_i,bt_i/a)_{\infty}}\\\times
\frac {(t_iq/c_i,t_iq/d_i,t_iq/e_i,
t_iq^{1-r+j}/f,t_iq^{2-j}/g, t_iq/h)_{\infty}}
{(c_it_i/a,d_it_i/a,e_it_i/a,
ft_iq^{r-j}/a,gt_iq^{j-1}/a, ht_i/a)_{\infty}}\;
d_qt_i\Bigg).
\end{multline*}
Now, to the integral inside the determinant we apply the
$q$-integral transformation \eqref{int109gl},
with the substitution $t\mapsto t_ix_i$, and the
replacements $a\mapsto ax_i^2$, $b\mapsto bx_i$,
$c\mapsto c_ix_i$, $d\mapsto d_ix_i$, $e\mapsto e_ix_i$,
$f\mapsto fq^{r-j}x_i$,
$g\mapsto gq^{j-1}x_i$, and $h\mapsto hx_i$. Thus we obtain
\begin{multline*}
\left(\frac ag\right)^{\binom r2}q^{-\binom r3}
\prod_{i=1}^r(aq^{2-r}/fg,fq^{2+r-2i}/g)_{i-1}^{-1}\\\times
\det_{1\le i,j\le r}\Bigg(\frac a{\lambda}
\frac{(b/ax_i,ax_iq/b,\lambda c_ix_i/a,
\lambda d_ix_i/a,\lambda e_ix_i/a)_\infty}
{(b/\lambda x_i,\lambda x_iq/b,c_ix_i,d_ix_i,e_ix_i)_\infty}\\\times
\frac{(bfq^{r-j}/\lambda,bgq^{j-1}/\lambda,bh/\lambda)_\infty}
{(bfq^{r-j}/a,bgq^{j-1}/a, bh/a)_\infty}\;
\int_{\lambda x_i}^b
\frac {(1-t_i^2/a)(t_iq/\lambda x_i,t_iq/b)_{\infty}}
{(t_ix_i,bt_i/\lambda)_{\infty}}\\\times
\frac {(at_iq/c_i\lambda,at_iq/d_i\lambda,at_iq/e_i\lambda,
t_iq^{1-r+j}/f,t_iq^{2-j}/g, t_iq/h)_{\infty}}
{(c_it_i/a,d_it_i/a,e_it_i/a,
ft_iq^{r-j}/\lambda,gt_iq^{j-1}/\lambda,ht_i/\lambda)_{\infty}}\;
d_qt_i\Bigg).
\end{multline*}
Now, by using linearity of the determinant with respect to rows and columns,
we take some factors out of the determinant and obtain
\begin{multline}\label{xgl}
\left(\frac a{f\lambda}\right)^{\binom r2}q^{-\binom r3}
\left(\frac a{\lambda}\right)^r
\prod_{i=1}^r(aq^{2-r}/fg,fq^{2+r-2i}/g)_{i-1}\\\times
\prod_{i=1}^r
\frac{(b/ax_i,ax_iq/b,\lambda c_ix_i/a,\lambda d_ix_i/a,\lambda e_ix_i/a,
bfq^{i-1}/\lambda,bgq^{i-1}/\lambda,bh/\lambda)_\infty}
{(b/\lambda x_i,\lambda x_iq/b,c_ix_i,d_ix_i,e_ix_i,
bfq^{i-1}/a,bgq^{i-1}/a, bh/a)_\infty}\\\times
\int_{\lambda x_1}^{b}\dots\int_{\lambda x_r}^{b}
\prod_{i=1}^r(1-t_i^2/\lambda)\frac{(qt_i/\lambda x_i,qt_i/b)_{\infty}}
{(t_ix_i,bt_i/\lambda)_{\infty}}\\\times
\prod_{i=1}^r
\frac{(aqt_i/c_i\lambda,aqt_i/d_i\lambda,
aqt_i/e_i\lambda,qt_i/f,qt_i/g,qt_i/h)_{\infty}\,t_i^{r-1}}
{(c_it_i/a,d_it_i/a,e_it_i/a,ft_i/\lambda,
gt_i/\lambda,ht_i/\lambda)_{\infty}
(t_iq^{1+r-j/f},\lambda q^{2-r}/gt_i)_{r-1}}\\\times
\det_{1\le i,j\le r}\Bigg(\frac{(ft_i/\lambda,f/t_i)_{r-j}}
{(q^{2-r}t_i/g,\lambda q^{2-r}/gt_i)_{r-j}}\Bigg)
\;d_qt_r\dots d_qt_1.
\end{multline}
The determinant can be evaluated by means of Lemma~\ref{lemdet1} with
$X_i\mapsto t_i$, $A\mapsto f/\lambda$, $B\mapsto q^{2-r}/g$, and
$C\mapsto \lambda$; specifically
\begin{multline*}
\det_{1\le i,j\le r}\Bigg(\frac{(ft_i/\lambda,f/t_i)_{r-j}}
{(q^{2-r}t_i/g,\lambda q^{2-r}/gt_i)_{r-j}}\Bigg)\\=
f^{\binom r2}q^{\binom r3}
\prod_{1\le i<j\le r}(t_i-t_j)(1-t_it_j/\lambda)
\prod_{i=1}^r\frac{(\lambda q^{2-r}/fg,fq^{2+r-2i}/g)_{i-1}}
{(q^{2-r}t_i/g,\lambda q^{2-r}/gt_i)_{r-1}}t_i^{1-r}.
\end{multline*}
Now note that $bh/\lambda=aq^{2-r}/fg$ and
$bh/a=\lambda q^{2-r}/fg$ whence we have
\begin{equation*}
\prod_{i=1}^r\frac{(bh/\lambda)_{\infty}(\lambda q^{2-r}/fg)_{i-1}}
{(bh/a)_{\infty}(aq^{2-r}/fg)_{i-1}}=\prod_{i=1}^r
\frac{(bhq^{i-1}/\lambda)_{\infty}}{(bhq^{i-1}/a)_{\infty}}.
\end{equation*}
Substituting these calculations into \eqref{xgl}
we arrive at the right-hand side of \eqref{cnint109gl}. 
\end{proof}

\begin{Remark}\label{remconv}
Note that in \eqref{cnint109gl}, the multiple $q$-integrals
converge absolutely since, when writing the left- and right-hand
sides as determinants of single $q$-integrals as in the above proof,
we have absolute convergence in each entry of the respective
determinants.
\end{Remark}

We can iterate the transformation in Theorem~\ref{cnint109}
to obtain a different transformation formula. We need to
permute some variables first, else we would get back the
original identity. In order to proceed, we specialize
the variables first and let $c_i=c$, $d_i=d$ and $e_i=e/x_i$,
for $i=1,\dots,r$. For convenience, we further interchange
the variables $e$ and $g$. In this first step, this gives the transformation
\begin{multline}\label{cngl}
\int_{ax_1}^{b}\dots\int_{ax_r}^{b}
\prod_{1\le i<j\le r}(t_i-t_j)(1-t_it_j/a)
\prod_{i=1}^r(1-t_i^2/a)\frac{(qt_i/ax_i,qt_i/b)_{\infty}}
{(t_ix_i,bt_i/a)_{\infty}}\\
\times\prod_{i=1}^r
\frac{(qt_i/c,qt_i/d,qt_i/e,qt_i/f,
qt_ix_i/g,qt_i/h)_{\infty}}
{(ct_i/a,dt_i/a,et_i/a,ft_i/a,
gt_i/ax_i,ht_i/a)_{\infty}}\;d_qt_r\dots d_qt_1\\
=\left(\frac a\lambda\right)^{\binom{r+1}2}
\prod_{i=1}^r\frac
{(b/ax_i,ax_iq/b,\lambda cx_i/a,\lambda dx_i/a,\lambda g/a)_\infty}
{(b/\lambda x_i,\lambda x_iq/b,cx_i,dx_i,g)_\infty}\\\times
\prod_{i=1}^r\frac
{(beq^{i-1}/\lambda,bfq^{i-1}/\lambda,bhq^{i-1}/\lambda)_\infty}
{(beq^{i-1}/a,bfq^{i-1}/a, bhq^{i-1}/a)_\infty}\\\times
\int_{\lambda x_1}^{b}\dots\int_{\lambda x_r}^{b}
\prod_{1\le i<j\le r}(t_i-t_j)(1-t_it_j/\lambda)
\prod_{i=1}^r(1-t_i^2/\lambda)\frac{(qt_i/\lambda x_i,qt_i/b)_{\infty}}
{(t_ix_i,bt_i/\lambda)_{\infty}}\\
\times\prod_{i=1}^r
\frac{(aqt_i/c\lambda,aqt_i/d\lambda,
aqt_ix_i/g\lambda,qt_i/e,qt_i/f,qt_i/h)_{\infty}}
{(ct_i/a,dt_i/a,gt_i/ax_i,et_i/\lambda,
ft_i/\lambda,ht_i/\lambda)_{\infty}}\;d_qt_r\dots d_qt_1,
\end{multline}
where $a^3q^{3-r}=bcdefgh$ and $\lambda=a^2q/cdg$.
Now, we iterate \eqref{cngl}, in the second step with
$a\mapsto a^2q/cdg$, $c\mapsto e$, $d\mapsto f$,
$e\mapsto aq/dg$, $f\mapsto aq/cg$, and $g\mapsto aq/cd$.
The result is the following.
\begin{Corollary}\label{cnint109c}
Let $a^3q^{3-r}=bcdefgh$. Then there holds
\begin{multline}\label{cnint109cgl}
\int_{ax_1}^{b}\dots\int_{ax_r}^{b}
\prod_{1\le i<j\le r}(t_i-t_j)(1-t_it_j/a)
\prod_{i=1}^r(1-t_i^2/a)\frac{(qt_i/ax_i,qt_i/b)_{\infty}}
{(t_ix_i,bt_i/a)_{\infty}}\\
\times\prod_{i=1}^r
\frac{(qt_i/c,qt_i/d,qt_i/e,qt_i/f,
qt_ix_i/g,qt_i/h)_{\infty}}
{(ct_i/a,dt_i/a,et_i/a,ft_i/a,
gt_i/ax_i,ht_i/a)_{\infty}}\;d_qt_r\dots d_qt_1\\
=\prod_{i=1}^r\frac
{(b/ax_i,ax_iq/b,ax_iq/cg,ax_iq/dg,ax_iq/eg,ax_iq/fg,
bhq^{r-1}/a)_\infty}
{(cx_i,dx_i,ex_i,fx_i,g,gq^{1-r}/hx_i,hq^rx_i/g)_\infty}\\\times
\prod_{i=1}^r\frac
{(gq^{1-i},aq^{2-i}/ch,aq^{2-i}/dh,aq^{2-i}/eh,
aq^{2-i}/fh)_\infty}
{(bcq^{i-1}/a,bdq^{i-1}/a,beq^{i-1}/a,bfq^{i-1}/a,
bhq^{i-1}/a)_\infty}\cdot
\left(\frac{agq^{1-r}}{bh}\right)^{\binom{r+1}2}
\\\times
\int_{bhq^{r-1}x_1/g}^{b}\dots\int_{bhq^{r-1}x_r/g}^{b}
\prod_{1\le i<j\le r}(t_i-t_j)(1-gq^{1-r}t_it_j/bh)\\\times
\prod_{i=1}^r(1-gq^{1-r}t_i^2/bh)
\frac{(gq^{2-r}t_i/bhx_i,qt_i/b,cgt_i/a,dgt_i/a)_{\infty}}
{(t_ix_i,gq^{1-r}t_i/h,at_iq^{2-r}/bch,
at_iq^{2-r}/bdh)_{\infty}}\\
\times\prod_{i=1}^r
\frac{(egt_i/a,fgt_i/a,
aq^{2-r}t_ix_i/bh,qt_i/h)_{\infty}}
{(at_iq^{2-r}/beh,
at_iq^{2-r}/bfh,gt_i/ax_i,gq^{1-r}t_i/b)_{\infty}}\;d_qt_r\dots d_qt_1.
\end{multline}
\end{Corollary}

\section{Multivariable nonterminating ${}_{10}\phi_9$
transformations}
\label{secm109}

We write out \eqref{mqint} explicitly for the integral on the left-hand side
of \eqref{cnint109gl}:
\begin{multline*}
\int_{ax_1}^{b}\dots\int_{ax_r}^{b}\prod_{1\le i<j\le r}(t_i-t_j)(1-t_it_j/a)
\prod_{i=1}^r(1-t_i^2/a)\frac{(qt_i/ax_i,qt_i/b)_{\infty}}
{(t_ix_i,bt_i/a)_{\infty}}\\
\times\prod_{i=1}^r
\frac{(qt_i/c_i,qt_i/d_i,qt_i/e_i,qt_i/f,
qt_i/g,qt_i/h)_{\infty}}
{(c_it_i/a,d_it_i/a,e_it_i/a,ft_i/a,
gt_i/a,ht_i/a)_{\infty}}\;d_qt_r\dots d_qt_1\\
=\sum_{S\subseteq\{1,2,\dots,r\}}
(-1)^{|S|}a^{|S|}b^{r-|S|}(1-q)^ra^{\binom{|S|}2}b^{\binom{r-|S|}2}
\prod_{i\in S}x_i\\\times
\sum_{k_1,\dots,k_r=0}^{\infty}
\underset{i,j\in S}{\prod_{1\le i<j\le r}}
(x_iq^{k_i}-x_jq^{k_j})(1-ax_ix_jq^{k_i+k_j})
\prod_{i\in S}(1-ax_i^2q^{2k_i})\\\times
\underset{i,j\notin S}{\prod_{1\le i<j\le r}}
(q^{k_i}-q^{k_j})(1-b^2q^{k_i+k_j}/a)
\prod_{i\notin S}(1-b^2q^{2k_i}/a)\\\times
\prod_{i\in S, j\notin S}(ax_iq^{k_i}-bq^{k_j})(1-bx_iq^{k_i+k_j})
(-1)^{\chi(i>j)}\\\times
\prod_{i\in S}
\frac{(q^{1+k_i},ax_iq^{1+k_i}/b,ax_iq^{1+k_i}/c_i,ax_iq^{1+k_i}/d_i,
ax_iq^{1+k_i}/e_i)_{\infty}}
{(ax_i^2q^{k_i},bx_iq^{k_i},c_ix_iq^{k_i},d_ix_iq^{k_i},
e_ix_iq^{k_i})_{\infty}}\\\times
\prod_{i\in S}
\frac{(ax_iq^{1+k_i}/f,ax_iq^{1+k_i}/g,ax_iq^{1+k_i}/h)_{\infty}}
{(fx_iq^{k_i},gx_iq^{k_i},hx_iq^{k_i})_{\infty}}
\prod_{i\notin S}
\frac{(bq^{1+k_i}/ax_i,q^{1+k_i})_{\infty}}
{(bx_iq^{k_i},b^2q^{k_i}/a)_{\infty}}\\\times
\prod_{i\notin S}
\frac{(bq^{1+k_i}/c_i,bq^{1+k_i}/d_i,
bq^{1+k_i}/e_i,bq^{1+k_i}/f,bq^{1+k_i}/g,bq^{1+k_i}/h)_{\infty}}
{(bc_iq^{k_i}/a,bd_iq^{k_i}/a,be_iq^{k_i}/a,
bfq^{k_i}/a,bgq^{k_i}/a,bhq^{k_i}/a)_{\infty}}\;\cdot q^{\sum_{i=1}^r k_i},
\end{multline*}
where $|S|$ denotes the number of elements of $S$, and $\chi$
is the truth function (which evaluates to one if the argument is true
and to zero otherwise). 
A similar expression is obtained for
the right-hand side of \eqref{cnint109gl}.
Now, if we divide both sides by
\begin{multline*}
(-1)^ra^{\binom{r+1}2}x_1\dots x_r(1-q)^r
\prod_{1\le i<j\le r}(x_i-x_j)(1-ax_ix_j)\\\times
\prod_{i=1}^r\frac{(q,ax_iq/b,ax_iq/c_i,ax_iq/d_i,ax_iq/e_i,
ax_iq/f,ax_iq/g,ax_iq/h)_{\infty}}
{(ax_i^2q,bx_i,c_ix_i,d_ix_i,e_ix_i,fx_i,gx_i,hx_i)_{\infty}}
\end{multline*}
and simplify, we obtain the following result which reduces to
\eqref{109ntgl} when $r=1$.

\begin{Corollary}[A $C_r$ nonterminating ${}_{10}\phi_9$ transformation]
\label{cnnt109}
Let 
$$a^3q^{3-r}=bc_id_ie_ix_ifgh$$
 and $\lambda=a^2q/c_id_ie_ix_i$
 for $i=1,\dots,r$. Then there holds
\begin{multline}\label{cnnt109gl}
\sum_{S\subseteq\{1,2,\dots,r\}}
\left(\frac ba\right)^{\binom{r-|S|}2}
\prod_{i\notin S}
\frac{(ax_i^2q,c_ix_i,d_ix_i,e_ix_i)_{\infty}}
{(ax_i/b,ax_iq/c_i,ax_iq/d_i,ax_iq/e_i)_{\infty}}\\\times
\prod_{i\notin S}
\frac{(fx_i,gx_i,hx_i,b/ax_i,bq/c_i,bq/d_i,bq/e_i,bq/f,bq/g,bq/h)_{\infty}}
{(ax_iq/f,ax_iq/g,ax_iq/h,b^2q/a,bc_i/a,bd_i/a,be_i/a,
bf/a,bg/a,bh/a)_{\infty}}\\\times
\sum_{k_1,\dots,k_r=0}^{\infty}
\underset{i,j\in S}{\prod_{1\le i<j\le r}}
\frac{(x_iq^{k_i}-x_jq^{k_j})(1-ax_ix_jq^{k_i+k_j})}
{(x_i-x_j)(1-ax_ix_j)}
\prod_{i\in S}\frac{(1-ax_i^2q^{2k_i})}{(1-ax_i^2)}\\\times
\underset{i,j\notin S}{\prod_{1\le i<j\le r}}
\frac{(q^{k_i}-q^{k_j})(1-b^2q^{k_i+k_j}/a)}
{(x_i-x_j)(1-ax_ix_j)}
\prod_{i\notin S}\frac{(1-b^2q^{2k_i}/a)}{(1-b^2/a)}\\\times
\prod_{i\in S, j\notin S}
\frac{(x_iq^{k_i}-bq^{k_j}/a)(1-bx_iq^{k_i+k_j})}
{(x_i-x_j)(1-ax_ix_j)}\\\times
\prod_{i\in S}\frac{(ax_i^2,bx_i,c_ix_i,d_ix_i,e_ix_i,fx_i,gx_i,hx_i)_{k_i}}
{(q,ax_iq/b,ax_iq/c_i,ax_iq/d_i,ax_iq/e_i,
ax_iq/f,ax_iq/g,ax_iq/h)_{k_i}}\\\times
\prod_{i\notin S}\frac{(b^2/a,bx_i,bc_i/a,bd_i/a,be_i/a,
bf/a,bg/a,bh/a)_{k_i}}
{(q,bq/ax_i,bq/c_i,bq/d_i,bq/e_i,bq/f,bq/g,bq/h)_{k_i}}\;
\cdot q^{\sum_{i=1}^r k_i}\\
=\prod_{i=1}^r\frac{(ax_i^2q,b/ax_i,
\lambda x_iq/f,\lambda x_iq/g,\lambda x_iq/h,
bfq^{i-1}/\lambda,bgq^{i-1}/\lambda,bhq^{i-1}/\lambda)_{\infty}}
{(\lambda x_i^2q,b/\lambda x_i,ax_iq/f,ax_iq/g,ax_iq/h,
bfq^{i-1}/a,bgq^{i-1}/a,bhq^{i-1}/a)_{\infty}}\\\times
\prod_{1\le i<j\le r}\frac{(1-\lambda x_ix_j)}{(1-ax_ix_j)}\;
\sum_{S\subseteq\{1,2,\dots,r\}}
\left(\frac b{\lambda}\right)^{\binom{r-|S|}2}\\\times
\prod_{i\notin S}
\frac{(\lambda x_i^2q,\lambda c_ix_i/a,\lambda d_ix_i/a,\lambda e_ix_i/a,
fx_i,gx_i,hx_i)_{\infty}}
{(\lambda x_i/b,ax_iq/c_i,ax_iq/d_i,ax_iq/e_i,
\lambda x_iq/f,\lambda x_iq/g,\lambda x_iq/h)_{\infty}}\\\times
\prod_{i\notin S}
\frac{(b/\lambda x_i,abq/c_i\lambda,abq/d_i\lambda,abq/e_i\lambda,
bq/f,bq/g,bq/h)_{\infty}}
{(b^2q/\lambda,bc_i/a,bd_i/a,be_i/a,
bf/\lambda,bg/\lambda,bh/\lambda)_{\infty}}\\\times
\sum_{k_1,\dots,k_r=0}^{\infty}
\underset{i,j\in S}{\prod_{1\le i<j\le r}}
\frac{(x_iq^{k_i}-x_jq^{k_j})(1-\lambda x_ix_jq^{k_i+k_j})}
{(x_i-x_j)(1-\lambda x_ix_j)}
\prod_{i\in S}\frac{(1-\lambda x_i^2q^{2k_i})}{(1-\lambda x_i^2)}\\\times
\underset{i,j\notin S}{\prod_{1\le i<j\le r}}
\frac{(q^{k_i}-q^{k_j})(1-b^2q^{k_i+k_j}/\lambda)}
{(x_i-x_j)(1-\lambda x_ix_j)}
\prod_{i\notin S}\frac{(1-b^2q^{2k_i}/\lambda)}{(1-b^2/\lambda)}\\\times
\prod_{i\in S, j\notin S}
\frac{(x_iq^{k_i}-bq^{k_j}/\lambda)(1-bx_iq^{k_i+k_j})}
{(x_i-x_j)(1-\lambda x_ix_j)}\\\times
\prod_{i\in S}\frac{(\lambda x_i^2,bx_i,\lambda c_ix_i/a,
\lambda d_ix_i/a,\lambda e_ix_i/a,fx_i,gx_i,hx_i)_{k_i}}
{(q,\lambda x_iq/b,ax_iq/c_i,ax_iq/d_i,ax_iq/e_i,
\lambda x_iq/f,\lambda x_iq/g,\lambda x_iq/h)_{k_i}}\\\times
\prod_{i\notin S}\frac{(b^2/\lambda,bx_i,bc_i/a,bd_i/a,be_i/a,
bf/\lambda,bg/\lambda,bh/\lambda)_{k_i}}
{(q,bq/\lambda x_i,abq/c_i\lambda,abq/d_i\lambda,abq/e_i\lambda,
bq/f,bq/g,bq/h)_{k_i}}\;
\cdot q^{\sum_{i=1}^r k_i}.
\end{multline}
\end{Corollary}

For the absolute convergence of the multiple series in \eqref{cnnt109gl},
see Remark~\ref{remconv}.

Similarly, we write out the transformation of sums resulting from
Corollary~\ref{cnint109c}. For $r=1$, the following Corollary reduces
to Bailey's transformation formula in \cite[Eq.~(8.1)]{bail10}
(cf.~\cite[Ex.~2.30]{grhyp}).

\begin{Corollary}[A $C_r$ nonterminating ${}_{10}\phi_9$ transformation]
\label{cnnt109c}
Let $a^3q^{3-r}=bcdefgh$. Then there holds
\begin{multline}\label{cnnt109cgl}
\sum_{S\subseteq\{1,2,\dots,r\}}
\left(\frac ba\right)^{\binom{r-|S|}2}
\prod_{i\notin S}
\frac{(ax_i^2q,cx_i,dx_i,ex_i,fx_i)_{\infty}}
{(ax_i/b,ax_iq/c,ax_iq/d,ax_iq/e,ax_iq/f)_{\infty}}\\\times
\prod_{i\notin S}
\frac{(g,hx_i,b/ax_i,bq/c,bq/d,bq/e,bq/f,bx_iq/g,bq/h)_{\infty}}
{(ax_i^2q/g,ax_iq/h,b^2q/a,bc/a,bd/a,be/a,bf/a,bg/ax_i,bh/a)_{\infty}}\\\times
\sum_{k_1,\dots,k_r=0}^{\infty}
\underset{i,j\in S}{\prod_{1\le i<j\le r}}
\frac{(x_iq^{k_i}-x_jq^{k_j})(1-ax_ix_jq^{k_i+k_j})}
{(x_i-x_j)(1-ax_ix_j)}
\prod_{i\in S}\frac{(1-ax_i^2q^{2k_i})}{(1-ax_i^2)}\\\times
\underset{i,j\notin S}{\prod_{1\le i<j\le r}}
\frac{(q^{k_i}-q^{k_j})(1-b^2q^{k_i+k_j}/a)}
{(x_i-x_j)(1-ax_ix_j)}
\prod_{i\notin S}\frac{(1-b^2q^{2k_i}/a)}{(1-b^2/a)}\\\times
\prod_{i\in S, j\notin S}
\frac{(x_iq^{k_i}-bq^{k_j}/a)(1-bx_iq^{k_i+k_j})}
{(x_i-x_j)(1-ax_ix_j)}\\\times
\prod_{i\in S}\frac{(ax_i^2,bx_i,cx_i,dx_i,ex_i,fx_i,g,hx_i)_{k_i}}
{(q,ax_iq/b,ax_iq/c,ax_iq/d,ax_iq/e,
ax_iq/f,ax_i^2q/g,ax_iq/h)_{k_i}}\\\times
\prod_{i\notin S}\frac{(b^2/a,bx_i,bc/a,bd/a,be/a,
bf/a,bg/ax_i,bh/a)_{k_i}}
{(q,bq/ax_i,bq/c,bq/d,bq/e,bq/f,bx_iq/g,bq/h)_{k_i}}\;
\cdot q^{\sum_{i=1}^r k_i}\\
=\prod_{i=1}^r\frac{(ax_i^2q,b/ax_i,bx_iq^r/g,
bchq^{r-1}x_i/a,bdhq^{r-1}x_i/a,behq^{r-1}x_i/a)_{\infty}}
{(bhx_i^2q^r/g,gq^{1-r}/hx_i,ax_iq/c,ax_iq/d,ax_iq/e,
ax_iq/f)_{\infty}}\\\times
\prod_{i=1}^r\frac{(bfhq^{r-1}x_i/a,
gq^{1-i},aq^{2-i}/ch,aq^{2-i}/dh,aq^{2-i}/eh,aq^{2-i}/fh)_{\infty}}
{(ax_iq/h,bcq^{i-1}/a,bdq^{i-1}/a,beq^{i-1}/a,
bfq^{i-1}/a,bhq^{i-1}/a)_{\infty}}\\
\times\prod_{1\le i<j\le r}\frac{(1-bhq^{r-1} x_ix_j/g)}{(1-ax_ix_j)}\;
\sum_{S\subseteq\{1,2,\dots,r\}}
\left(\frac{gq^{1-r}}{h}\right)^{\binom{r-|S|}2}\\\times
\prod_{i\notin S}
\frac{(bh x_i^2q^r/g,ax_iq/cg,ax_iq/dg)_{\infty}}
{(hx_iq^{r-1}/g,bchq^{r-1}x_i/a,bdhq^{r-1}x_i/a)_{\infty}}\\\times
\prod_{i\notin S}
\frac{(ax_iq/eg,ax_iq/fg,bhq^{r-1}/a,hx_i,gq^{1-r}/hx_i)_{\infty}}
{(behq^{r-1}x_i/a,bfhq^{r-1}x_i/a,ax_i^2q/g,
bx_iq^r/g,bgq^{2-r}/h)_{\infty}}\\\times
\prod_{i\notin S}
\frac{(bcg/a,bdg/a,beg/a,bfg/a,ax_iq^{2-r}/h,bq/h)_{\infty}}
{(aq^{2-r}/ch,aq^{2-r}/dh,aq^{2-r}/eh,aq^{2-r}/fh,
bg/ax_i,gq^{1-r})_{\infty}}\\\times
\sum_{k_1,\dots,k_r=0}^{\infty}
\underset{i,j\in S}{\prod_{1\le i<j\le r}}
\frac{(x_iq^{k_i}-x_jq^{k_j})(1-bh x_ix_jq^{r-1+k_i+k_j}/g)}
{(x_i-x_j)(1-bh x_ix_jq^{r-1}/g)}\\\times
\underset{i,j\notin S}{\prod_{1\le i<j\le r}}
\frac{(q^{k_i}-q^{k_j})(1-bgq^{1-r+k_i+k_j}/h)}
{(x_i-x_j)(1-bh x_ix_jq^{r-1}/g)}
\prod_{i\notin S}\frac{(1-bgq^{1-r+2k_i}/h)}{(1-bgq^{1-r}/h)}\\\times
\prod_{i\in S, j\notin S}
\frac{(x_iq^{k_i}-gq^{1-r+k_j}/h)(1-bx_iq^{k_i+k_j})}
{(x_i-x_j)(1-bhx_ix_jq^{1-r}/g)}
\prod_{i\in S}\frac{(1-bhx_i^2q^{r-1+2k_i}/g)}{(1-bh x_i^2q^{r-1}/g)}\\\times
\prod_{i\in S}\frac{(bh x_i^2q^{r-1}/g,bx_i,ax_iq/cg,
ax_iq/dg,ax_iq/eg,ax_iq/fg)_{k_i}}
{(q,hx_iq^r/b,bchq^{r-1}x_i/a,bdhq^{r-1}x_i/a,behq^{r-1}x_i/a,
bfhq^{r-1}x_i/a)_{k_i}}\\\times
\prod_{i\in S}\frac{(bhq^{r-1}/a,hx_i)_{k_i}}
{(ax_i^2q/g,bx_iq^r/g)_{k_i}}\;
\prod_{i\notin S}\frac{(bgq^{1-r}/h,bx_i)_{k_i}}
{(q,gq^{2-r}/hx_i)_{k_i}}\\\times
\prod_{i\notin S}\frac{(aq^{2-r}/ch,aq^{2-r}/dh,
aq^{2-r}/eh,aq^{2-r}/fh,bg/ax_i,gq^{1-r})_{k_i}}
{(bcg/a,bdg/a,beg/a,bfg/a,ax_iq^{2-r}/h,bq/h)_{k_i}}\;
\cdot q^{\sum_{i=1}^r k_i}.
\end{multline}
\end{Corollary}

\section{Terminating $C_r$ ${}_{10}\phi_9$ transformations}\label{sect109}

We will now specialize  Corollary~\ref{cnnt109} to
obtain terminating multivariable ${}_{10}\phi_9$
transformations. We accomplish this by 
multiplying both sides of the identity with
$\prod_{i=1}^r(bd_i/a)_{\infty}$ and then letting $bd_i/a=q^{-n_i}$. Then
only the terms corresponding to $S=\emptyset$ are non-zero.
After a change of variables this gives the following transformation:
\begin{Corollary}[A $C_r$ terminating ${}_{10}\phi_9$ transformation]
\label{cnt109}
Let $a^3q^{3-r+n_i}=bcde_if_ig_i$, for $i=1,\dots,r$, and
$\lambda=a^2q^{2-r}/bcd$. Then there holds
\begin{multline}\label{cnt109gl}
\underset{i=1,\dots,r}{\sum_{0\le k_i\le n_i}}
\prod_{1\le i<j\le r}(q^{k_i}-q^{k_j})
  (1-aq^{k_i+k_j})\\
\times \prod_{i=1}^r\frac{(1-aq^{2k_i})(a,b,c,d,e_i,f_i,g_i,q^{-n_i})_{k_i}}
{(1-a)(q,aq/b,aq/c,aq/d,aq/e_i,aq/f_i,aq/g_i,aq^{1+n_i})_{k_i}}\,
q^{k_i}\\
=\left(\frac\lambda a\right)^{\binom r2}
\prod_{i=1}^{r-1}\frac{(b,c,d)_i}
{(b\lambda/a,c\lambda/a,d\lambda/a)_i}\prod_{i=1}^r\frac
{(aq,\lambda q/e_i, \lambda q/f_i, aq/e_if_i)_{n_i}}
{(\lambda q,aq/e_i, aq/f_i,\lambda q/e_if_i)_{n_i}}\\
\times \underset{i=1,\dots,r}{\sum_{0\le k_i\le n_i}}
\prod_{1\le i<j\le r}(q^{k_i}-q^{k_j})(1-\lambda q^{k_i+k_j})\\
\times \prod_{i=1}^r\frac{(1-\lambda q^{2k_i})
(\lambda ,\lambda b/a,\lambda c/a,\lambda d/a,e_i, f_i,g_i,q^{-n_i})_{k_i}}
{(1-\lambda)(q,aq/b,aq/c,aq/d,\lambda q/e_i,\lambda q/f_i,\lambda
  q/g_i,\lambda q^{1+n_i})_{k_i}}\,q^{k_i}.
\end{multline}
\end{Corollary}
Note that in \eqref{cnt109gl} only terms with all $k_i$ distinct are non-zero.
 Corollary~\ref{cnt109} reduces to Bailey's~\cite{bail10}
transformation formula for $r=1$ (cf.~\cite[Eq.~(2.9.1)]{grhyp}).

Similarly to the derivation of Corollary~\ref{cnnt109c} from
Corollary~\ref{cnnt109}, we can obtain another transformation
from \eqref{cnt109gl} by iteration. For this, we first
specialize the variables letting $f_i=f$, $g_i=g$,
for $i=1,\dots,r$, and write $e_i=eq^{n_i}$ for convenience.
We then iterate the relation with
$a\mapsto a^2q^{2-r}/bcd$, $b\mapsto aq^{2-r}/cd$, $c\mapsto g$, $d\mapsto f$,
$f\mapsto aq^{2-r}/bc$, and $g\mapsto aq^{2-r}/bd$.
We obtain the following result which reduces to 
\cite[Ex.~2.19]{grhyp} (with $e$ replaced by $eq^n$)
for $r=1$.
\begin{Corollary}
[A $C_r$ terminating ${}_{10}\phi_9$ transformation]\label{cnt109n}
Let $bcdefg=a^3q^{3-r}$. Then there holds
\begin{multline}\label{cnt109ngl}
\underset{i=1,\dots,r}{\sum_{0\le k_i\le n_i}}
\prod_{1\le i<j\le r}(q^{k_i}-q^{k_j})(1-aq^{k_i+k_j})\\
\times \prod_{i=1}^r\frac{(1-aq^{2k_i})(a,b,c,d,eq^{n_i},f,g,q^{-n_i})_{k_i}}
{(1-a)(q,aq/b,aq/c,aq/d,aq^{1-n_i}/e,aq/f,aq/g,aq^{1+n_i})_{k_i}}\,
q^{\sum_{i=1}^rk_i}\\
=\left(\frac{a^2q^{4-2r}}{b^2cdfg}\right)^{\binom r2}
\prod_{i=1}^{r-1}\frac{(b,c,d,f,g)_i}
{(aq^{2-r}/bc,aq^{2-r}/bd,aq^{2-r}/bf,aq^{2-r}/bg,a^2q^{4-2r}/bcdfg)_i}\\\times
\prod_{i=1}^r\frac
{(aq,aq^{1-n_i}/ce,aq^{1-n_i}/de,aq^{1-n_i}/ef,aq^{1-n_i}/eg,bq^{r-1})_{n_i}}
{(aq/c,aq/d,aq^{1-n_i}/e,aq/f,aq/g,bq^{r-1-n_i}/e)_{n_i}}(eq^{n_i})^{n_i}\\
\times \underset{i=1,\dots,r}{\sum_{0\le k_i\le n_i}}
\prod_{1\le i<j\le r}(q^{k_i}-q^{k_j})(1-eq^{1-r+k_i+k_j}/b)\\
\times \prod_{i=1}^r\frac{(1-eq^{1-r+2k_i}/b)
(eq^{1-r}/b,eq^{1-r}/a,aq^{2-r}/bf,aq^{2-r}/bg)_{k_i}}
{(1-eq^{1-r}/b)(q,aq/b,ef/a,eg/a)_{k_i}}\\
\times \prod_{i=1}^r\frac{(eq^{n_i},aq^{2-r}/bc,aq^{2-r}/bd,q^{-n_i})_{k_i}}
{(q^{2-r-n_i}/b,ec/a,ed/a,eq^{2-r+n_i}/b)_{k_i}}\,
q^{\sum_{i=1}^rk_i}.
\end{multline}
\end{Corollary}
Again, only terms with all $k_i$ distinct are non-zero.

Finally, we give another multivariable extension of \cite[Ex.~2.19]{grhyp}.
We first let $g=q^{-N}$ in Corollary~\ref{cnnt109c} and then do the
simultaneous replacements $b\mapsto e$, $e\mapsto g$, and $h\mapsto b$.

\begin{Corollary}[A $C_r$ terminating ${}_{10}\phi_9$ transformation]
\label{cnt109i2}
Let $a^3q^{3-r+N}=bcdefg$. Then there holds
\begin{multline}\label{cnt109i2gl}
\sum_{k_1,\dots,k_r=0}^N
\prod_{1\le i<j\le r}
\frac{(x_iq^{k_i}-x_jq^{k_j})(1-ax_ix_jq^{k_i+k_j})}
{(x_i-x_j)(1-ax_ix_j)}
\prod_{i=1}^r\frac{(1-ax_i^2q^{2k_i})}{(1-ax_i^2)}\,q^{k_i}\\
\times
\prod_{i=1}^r\frac{(ax_i^2,bx_i,cx_i,dx_i,ex_i,fx_i,gx_i,q^{-N})_{k_i}}
{(q,ax_iq/b,ax_iq/c,ax_iq/d,ax_iq/e,ax_iq/e,ax_iq/f,ax_iq/g,
ax_i^2q^{1+N})_{k_i}}\\
=
q^{\binom r 3}\left(\frac{q^{N+1}}{b}\right)^{\binom r 2}\frac{1}
{(b/e)_{N+r-1}^r}
\prod_{1\le i<j\le r}
\frac{1}{(x_i-x_j)(1-ax_ix_j)}\\
\times\prod_{i=1}^{r}\frac {(bx_i)^N
(ax_i^2q)_N(bx_i)_{N+r-1}(q^{2-r}ax_i/b)_{r-1}}
{(q^{N+r+1-i},aq^{2-r}/be,aq^{2-r}/bc,aq^{2-r}/bd,
aq^{2-r}/bf,aq^{2-r}/bg)_{i-1}}\\
\times\prod_{i=1}^r\frac{(q^{2-i}a/ce,q^{2-i}a/de,q^{2-i}a/ef,
 q^{2-i}a/eg)_{N+i-1}}
{(ax_iq/c,ax_iq/d,ax_iq/e,ax_iq/f,ax_iq/g)_N}\\
\times
\sum_{k_1,\dots,k_r=0}^{N+r-1}\prod_{1\le i<j\le r}
(q^{k_i}-q^{k_j})(1-eq^{1-r-N+k_i+k_j}/b)\\
\times
\prod_{i=1}^r
\frac{(1-eq^{1-r-N+2k_i}/b)(eq^{1-r-N}/b,ex_i,eq^{-N}/ax_i,q^{1-r-N})_{k_i}}
{(1-eq^{1-r-N}/b)(q,q^{2-r-N}/bx_i,ax_iq^{2-r}/b,eq/b)_{k_i}}\\
\times\prod_{i=1}^r\frac{(aq^{2-r}/bc,aq^{2-r}/bd,aq^{2-r}/bf,
aq^{2-r}/bg)_{k_i}}
{(ceq^{-N}/a,deq^{-N}/a,efq^{-N}/a,egq^{-N}/a)_{k_i}}\,q^{k_i}.
\end{multline}
\end{Corollary}

Note that the sum on the right lives on a larger hypercube, though
this is to some extent ``compensated'' by the fact that, on the right
but not on the left, only terms with all $k_i$ distinct are non-zero.

\section{$A_{r-1}$ ${}_1\psi_1$ summations}
\label{secmult}

Previously, one of us~\cite[Theorems~2.1 and 3.3]{schlhypdet}
utilized determinant evaluations to derive (among other summations)
several $A_{r-1}$ extensions of Ramanujan's ${}_1\psi_1$ summation.
Here we provide more multivariable extensions of
Ramanujan's ${}_1\psi_1$ summation by the same method.
The difference is that, like for most sums in Section
\ref{sect109}, the terms of 
our bilateral series are now zero unless all summation indices
$k_i$ are different.

We note that we were not able to
give new $C_r$ extensions of Bailey's ${}_6\psi_6$ summation
(cf.\ \cite[Eq.~(II.33)]{grhyp}) by the same method.
All such obtained ${}_6\psi_6$ summations would be in fact
just special cases of an identity given earlier
by one of us, namely the identity obtained from
equating Eq.~(3.4) with Eq.~(3.5) from \cite{schlhypdet}.

The following result reduces to Ramanujan's ${}_1\psi_1$ summation
(cf.\ \cite[Eq.~(5.2.1)]{grhyp}) for $r=1$.

\begin{Proposition}[An $A_{r-1}$ ${}_1\psi_1$ summation]\label{an11}
We have
\begin{multline}\label{an11gl}
\sum_{k_1,\dots,k_r=-\infty}^{\infty}
\prod_{1\le i<j\le r}(q^{k_i}-q^{k_j})
\prod_{i=1}^r\frac{(a_i)_{k_i}}{(b_i)_{k_i}}z_i^{k_i}\\
=
q^{\binom r2}
\det_{1\le i,j\le r}\left(b_i^{j-r}\frac{(z_i)_{r-j}}
{(a_iz_iq/b_i)_{r-j}}\right)
\prod_{i=1}^r\frac{(q,a_iz_i,q/a_iz_i,b_i/a_i)_{\infty}}
{(b_i,z_i,b_i/a_iz_i,q/a_i)_{\infty}},
\end{multline}
provided $|b_iq^{1-r}/a_i|<|z_i|<1$, for $i=1,\dots,r$.
\end{Proposition}

\begin{proof}
We have
\begin{equation*}
\prod_{1\le i<j\le r}(q^{k_i}-q^{k_j})=
\det_{1\le i,j\le r}\left(\big(q^{k_i}\big)^{r-j}\right),
\end{equation*}
due to the classical Vandermonde determinant evaluation.
Hence we may write the left-hand side of \eqref{an11gl} as
\begin{equation*}
\det_{1\le i,j\le r}\left(\sum _{k_i=-\infty}^{\infty}
\frac{(a_i)_{k_i}}{(b_i)_{k_i}}\left(z_iq^{r-j}\right)^{k_i}\right).
\end{equation*}
Now, to the sum inside the determinant we apply
Ramanujan's $_1\psi_1$ summation \eqref{11gl},
with $a\mapsto a_i$, $b\mapsto b_i$, and $z\mapsto z_iq^{r-j}$.
Thus we obtain
\begin{equation*}
\det_{1\le s,t\le r}\left(
\frac{(q,a_iz_iq^{r-j},q^{1+j-r}/a_iz_i,b_i/a_i)_{\infty}}
{(b_i,z_iq^{r-j},b_iq^{j-r}/a_iz_i,q/a_i)_{\infty}}
\right).
\end{equation*}
Now, by using linearity of the determinant with respect to rows, we take
some factors out of the determinant and readily obtain
the right-hand side of \eqref{an11gl}.
\end{proof}

We now consider different
specializations of the parameters $a_i$, $b_i$, and $z_i$, for $i=1,\dots,r$,
for which the determinant in Proposition~\ref{an11}
can be reduced to a product by means of
Lemma~\ref{lemdet1}.

A simple choice would be $b_i=b$ and $z_i=z$, for $i=1,\dots,r$. In this
case we recover the special case of Eq.~(2.2) of \cite[Th.~2.1]{schlhypdet}
where in the latter identity both sides are multiplied by
$\prod_{1\le i<j\le r}(1-x_i/x_j)$ and then the specializations
$x_i=1$, for $i=1,\dots,r$, are being made.
Another simple choice would be $a_i=a$ and $b_i=b$, for $i=1,\dots,r$.
In this case we similarly recover a special case of Eq.~(2.3) of
\cite[Th.~2.1]{schlhypdet}.

We also recover a previous result if in Proposition~\ref{an11} we choose the
specializations $a_i=a/x_i$, $b_i=b$, and $z_i=zx_i$, for $i=1,\dots,r$.
This yields a special case of Eq.~(2.10) of \cite{schlhypdet}.
Similarly, if we choose $a_i=a$, $b_i=bx_i$, and $z_i=z$, for $i=1,\dots,r$,
the determinant can be evaluated by \eqref{lemdet1a} and we obtain a
special case of Eq.~(2.13) of \cite{schlhypdet}.

Nevertheless, we are able to give two different multivariable
${}_1\psi_1$ sums resulting from special cases of
Proposition~\ref{an11} which are not covered by previous results.

First, we choose $a_i=ax_i$, $b_i=bx_i$, and $z_i=z$, for $i=1,\dots,r$.
In this case we obtain the following result.

\begin{Corollary}[An $A_{r-1}$ ${}_1\psi_1$ summation]\label{an11b}
We have
\begin{multline}
\sum_{k_1,\dots,k_r=-\infty}^{\infty}
\prod_{1\le i<j\le r}(q^{k_i}-q^{k_j})
\prod_{i=1}^r\frac{(ax_i)_{k_i}}{(bx_i)_{k_i}}z^{k_i}\\
=
(az)^{-\binom r2}q^{-\binom r3}\prod_{1\le i<j\le r}(1/x_j-1/x_i)
\prod_{i=1}^r
\frac{(q,azx_i,q/azx_i,b/a)_{\infty}}
{(bx_i,zq^{i-1},bq^{1-i}/az,q/ax_i)_{\infty}},
\end{multline}
provided $|bq^{1-r}/a|<|z|<1$.
\end{Corollary}

Next, we choose $a_i=a/x_i$, $b_i=bx_i$, and $z_i=zx_i$, for $i=1,\dots,r$.
We now require a particular limit case of the determinant evaluation
in Lemma~\ref{lemdet1}, which can be generally stated as
\begin{equation}\label{lemdet1a}
\det_{1\le i,j\le r}\left(X_i^{r-j}
\frac{(A/X_i)_{r-j}}{(BX_i)_{r-j}}\right)
=
\prod_{1\le i<j\le r}(X_i-X_j)\prod_{i=1}^r
\frac{(ABq^{2r-2i})_{i-1}}{(BX_i)_{r-1}}.
\end{equation}
The result is the following.

\begin{Corollary}[An $A_{r-1}$ ${}_1\psi_1$ summation]\label{an11e}
We have
\begin{multline}
\sum_{k_1,\dots,k_r=-\infty}^{\infty}
\prod_{1\le i<j\le r}(q^{k_i}-q^{k_j})
\prod_{i=1}^r\frac{(a/x_i)_{k_i}}{(bx_i)_{k_i}}\big(zx_i\big)^{k_i}\\
=
a^{-\binom r2}
\prod_{1\le i<j\le r}(x_i-x_j)\prod_{i=1}^r
(bq^{1-2r+i}/az^2)_{i-1}
\frac{(q,az,q/az,bx_i^2/a)_{\infty}}
{(bx_i,zx_i,bx_iq^{2-r}/az,x_iq/a)_{\infty}},
\end{multline}
provided $|bx_i^2q^{1-r}/a|<|zx_i|<1$ for $i=1,\dots,r$.
\end{Corollary}

\section{Specializations}

We will now
 point out some interesting special cases and consequences
of our new $C_r$ ${}_{10}\phi_9$ transformations.

\subsection{A Watson transformation}
In Corollary~\ref{cnt109}, we first remove the dependency
of the parameters by replacing $g_i$ by
$a^3q^{3-r+n_i}/bcde_if_i$. If we now let $d\to\infty$ and
relabel $f_i\mapsto d_i$, for $i=1,\dots,r$, we obtain the
following multivariable generalization of Watson's transformation
(cf.\ \cite[Eq.~(2.5.1)]{grhyp} ):

\begin{Corollary}[A multivariable Watson transformation]\label{mwatson}
We have
\begin{multline}\label{mwatsongl}
\underset{i=1,\dots,r}{\sum_{0\le k_i\le n_i}}
\prod_{1\le i<j\le r}(q^{k_i}-q^{k_j})
  (1-aq^{k_i+k_j})\\
\times \prod_{i=1}^r\frac{(1-aq^{2k_i})(a,b,c,d_i,e_i,q^{-n_i})_{k_i}}
{(1-a)(q,aq/b,aq/c,aq/d_i,aq/e_i,aq^{1+n_i})_{k_i}}
\left(\frac{a^2q^{3-r+n_i}}{bcd_ie_i}\right)^{k_i}\\
=\left(-\frac a{bc}\right)^{\binom r2}q^{-2\binom r3}
\prod_{i=1}^{r-1}\frac{(b,c)_i}
{(aq^{2-r}/bc)_i}\prod_{i=1}^r\frac
{(aq, aq/d_ie_i)_{n_i}}
{(aq/d_i, aq/e_i)_{n_i}}\\
\times \underset{i=1,\dots,r}{\sum_{0\le k_i\le n_i}}
\prod_{1\le i<j\le r}(q^{k_i}-q^{k_j})
\prod_{i=1}^r\frac{(aq^{2-r}/bc,d_i,e_i,q^{-n_i})_{k_i}}
{(q,aq/b,aq/c,d_ie_iq^{-n_i}/a)_{k_i}}\,q^{k_i}.
\end{multline}
\end{Corollary}

In the special case when   $bc=aq$,
the right-hand side can be written as a multiple of a determinant;
cf.\  the proof of Corollary~\ref{cnt87}.  This gives the $C_r$
terminating ${}_6\phi_5$ summation stated as Corollary \ref{cnt65}
below. Similarly, the case $b=q^{1-r}$ gives the $A_{r-1}$ terminating
${}_3\phi_2$ summation in Corollary \ref{an32} below.
As a matter of fact, we can even extract a $C_r$ Jackson summation
from our Watson transformation in Corollary~\ref{mwatson}.
We specialize the parameters so that the series on the left-hand side
becomes balanced, i.e., we set $a^3q^{2-r+n_i}=bcd_ie_i$, for $i=1,\dots,r$.
Note that in this case the multivariable ${}_4\phi_3$ on the right-hand
side of \eqref{mwatsongl} reduces to a multivariable ${}_3\phi_2$
which can be transformed into a multiple of a
determinant by Corollary \ref{an32}.
(Here we can assume that Corollary \ref{an32} is already known --
after all, we just pointed out that it follows from
Corollary~\ref{mwatson}.)
However, the result we would obtain by this procedure is only
a special case of Corollary~\ref{cnt87} below.

We also remark that on the right-hand side, 
the ``type $C$'' Vandermonde determinant 
\begin{equation}\label{facdiff}
\prod_{1\le i<j\le r}(q^{k_i}-q^{k_j})(1-aq^{k_i+k_j})
\prod_{i=1}^r(1-aq^{2k_i})
\end{equation}
 got reduced to the classical ``type $A$''
Vandermonde determinant. Moreover, the very-well-poised
condition of the parameters got lost, while the series remains
 balanced. The left-hand side retains the factor \eqref{facdiff} and
 is very-well-poised but not balanced.
 Dealing here with $r$-dimensional series, 
we simply mean by these terms
that the respective series are very-well-poised and/or balanced when $r=1$.

\subsection{An $\bmar$ Sears transformation}

To obtain a multivariable extension of Sears' transformation
(cf.\ \cite[Eq.~(3.2.1)]{grhyp}) from Corollary~\ref{cnt109},
replace $b$ by $aq/b$ and $e_i$ by $aq/e_i$, for $i=1,\dots,r$,
and then take the limit $a\to 0$. After relabeling of parameters,
$b\mapsto d$, $d\mapsto b$, $f_i\mapsto a_i$, for $i=1,\dots,r$,
we obtain

\begin{Corollary}[An $A_{r-1}$ Sears transformation]\label{msears}
We have
\begin{multline}\label{msearsgl}
\underset{i=1,\dots,r}{\sum_{0\le k_i\le n_i}}
\prod_{1\le i<j\le r}(q^{k_i}-q^{k_j})
\prod_{i=1}^r\frac{(a_i,b,c,q^{-n_i})_{k_i}}
{(q,d,e_i,a_ibcq^{r-n_i}/de_i)_{k_i}}q^{k_i}\\
=\left(\frac{dq^{1-r}}{bc}\right)^{\binom r2}
\prod_{i=1}^{r-1}\frac{(b,c)_i}
{(dq^{1-r}/b,dq^{1-r}/c)_i}\prod_{i=1}^r\frac
{(de_iq^{-r}/bc, e_i/a_i)_{n_i}}
{(e_i,de_iq^{-r}/a_ibc)_{n_i}}\\\times 
\underset{i=1,\dots,r}{\sum_{0\le k_i\le n_i}}
\prod_{1\le i<j\le r}(q^{k_i}-q^{k_j})
\prod_{i=1}^r\frac{(a_i,dq^{1-r}/b,dq^{1-r}/c,q^{-n_i})_{k_i}}
{(q,d,de_iq^{-r}/bc,a_iq^{1-n_i}/e_i)_{k_i}}q^{k_i}.
\end{multline}
\end{Corollary}

When $d=b$, the right-hand side can be expressed as a determinant.
This gives a special case of Corollary \ref{an32} below.

\subsection{$\bmcr$ Jackson summations}\label{subsect87}
We start with the
following identity, which reduces to Jackson's summation formula \eqref{87gl}
when $r=1$.

\begin{Corollary}[A $C_r$ Jackson summation]\label{cnt87}
Let $bc_id_ie_i=a^2q^{2-r+n_i}$ for $i=1,\dots,r$. Then 
\begin{multline}\label{cnt87gl}
\underset{i=1,\dots,r}{\sum_{0\le k_i\le n_i}}
\prod_{1\le i<j\le r}(q^{k_i}-q^{k_j})
  (1-aq^{k_i+k_j})\\
\times \prod_{i=1}^r\frac{(1-aq^{2k_i})(a,b,c_i,d_i,e_i,q^{-n_i})_{k_i}}
{(1-a)(q,aq/b,aq/c_i,aq/d_i,aq/e_i,aq^{1+n_i})_{k_i}}\,
q^{k_i}\\
=(-b)^{-\binom r2}q^{-2\binom r 3}
\prod_{i=1}^r\frac{(aq^{2-r}/b)_{r-1}(b)_{i-1}}{(aq^{2+r-2i}/b)_{i-1}}
\frac{(aq,aq/c_id_i,aq/c_ie_i,aq/d_ie_i)_{n_i}}
{(aq/c_i,aq/d_i,aq/e_i,aq/c_id_ie_i)_{n_i}}\\
\times\det_{1\le i,j\le r}\left(
\frac{(c_i,d_i,e_i,q^{-n_i})_{r-j}}
{(aq^{2-r}/bc_i,aq^{2-r}/bd_i,aq^{2-r}/be_i,aq^{2-r+n_i}/b)_{r-j}}\right).
\end{multline}
\end{Corollary}

\begin{proof}
To obtain the left-hand side, we 
 let  $cd=aq$ in Corollary \ref{cnt109} and make the change of
variables  $(e_i,f_i,g_i)\mapsto(c_i,d_i,e_i)$.
 On the right-hand side we observe that, since
 $\lambda b/a=q^{1-r}$,  only terms with $k_i\leq r-1$ are
 non-zero. Since we may also assume that the
 $k_i$ are  distinct, $(k_1,\dots k_r)$ must be a permutation of
$(0,\dots,r-1)$. This allows us to pull out some factors from the
 sum, for instance
$$\prod_{1\le i<j\le r}(q^{k_i}-q^{k_j})=\prod_{0\le i<j\le
  r-1}(q^{i}-q^{j})\operatorname{sgn}(k), $$
where $\operatorname{sgn}$ denotes the sign of the permutation.
After some cancellation, we obtain the expression
\begin{multline}
\left(\frac{\lambda q} a\right)^{\binom r2}
\prod_{i=1}^{r-1}\frac{(1-\lambda q^{2i})(b,\lambda)_i}
{(1-\lambda)(q,aq/b)_i}\prod_{i=1}^r\frac
{(aq,\lambda q/c_i, \lambda q/d_i, aq/c_id_i)_{n_i}}
{(\lambda q,aq/c_i, aq/d_i,\lambda q/c_id_i)_{n_i}}\\
\times \prod_{0\le i<j\le
  r-1}(q^{i}-q^{j})(1-\lambda q^{i+j})
\sum_{k}\operatorname{sgn}(k)
 \prod_{i=1}^r\frac{
(c_i,d_i,e_i,q^{-n_i})_{k_i}}
{(\lambda q/c_i,\lambda q/d_i,\lambda q/e_i,
\lambda q^{1+n_i})_{k_i}}
\end{multline}
for the right-hand side. Up to a factor $(-1)^{\binom r 2}$, obtained
from inverting the order of the columns, the 
 sum in $k$ equals the determinant in \eqref{cnt87gl}.
After some manipulations of  $q$-shifted factorials, including the
easily verified identities
$$\prod_{0\le i<j\le r-1}(q^{i}-q^{j})\prod_{i=1}^{r-1}\frac{1}
{(q)_{i}}=q^{\binom r 3} $$
and (recall that  $\lambda=aq^{1-r}/b$)
$$\prod_{0\le i<j\le r-1}(1-\lambda
 q^{i+j})\prod_{i=1}^{r-1}\frac{(1-\lambda q^{2i})(\lambda)_{i}}
{(1-\lambda)(aq/b)_{i}}=(aq^{2-r}/b)_{r-1}^r\prod_{i=1}^{r-1}\frac{1}
{(aq^{r-2i}/b)_i}, $$
 we arrive at
 \eqref{cnt87gl}. 
\end{proof}

Next we consider two cases of
Corollary \ref{cnt87}  when the determinant
on the right-hand side is computed by Lemma~\ref{lemdet1}. 
For the first case we put
$c_i=c$, $d_i=d$ and write $e_i=eq^{n_i}$.

\begin{Corollary}[A $C_r$ Jackson summation]\label{cnt87a}
If $bcde=a^2q^{2-r}$, then 
\begin{multline}
\underset{i=1,\dots,r}{\sum_{0\le k_i\le n_i}}
\prod_{1\le i<j\le r}(q^{k_i}-q^{k_j})
  (1-aq^{k_i+k_j})\\
\times \prod_{i=1}^r\frac{(1-aq^{2k_i})(a,b,c,d,eq^{n_i},q^{-n_i})_{k_i}}
{(1-a)(q,aq/b,aq/c,aq/d,aq^{1-n_i}/e,aq^{1+n_i})_{k_i}}\,
q^{k_i}\\
=q^{-\binom r 3}\left(\frac e a\right)^{\binom r2}
\prod_{i=1}^{r-1}\frac{(b,c,d)_i}{(aq^{2-r}/bc,aq^{2-r}/bd,aq^{2-r}/cd)_i}\\
\times\prod_{1\le i<j\le r}(q^{n_i}-q^{n_j})
  (1-eq^{n_i+n_j})
\prod_{i=1}^r
\frac{(aq,aq^{2-r}/bc,aq^{2-r}/bd,aq^{2-r}/cd)_{n_i}}
{(aq/b,aq/c,aq/d,aq^{2-r}/bcd)_{n_i}}.
\end{multline}
\end{Corollary}

Next we give the case of Corollary \ref{cnt87} when $c_i=c$ and $n_i=N$ for all
$i$, and we write $d_i=dx_i$, $e_i=e/x_i$. We have used some
 manipulations to write the result so that the symmetry
between $b$ and $c$ is exhibited. An alternative way to obtain this
identity,  which gives it in the form we want immediately, is to put
$bh/a=q^{-N}$ and $d_ie_i=aq$ in Corollary \ref{cnnt109}.
 
\begin{Corollary}[A $C_r$ Jackson summation]\label{cnt87b}
If $bcde=a^2q^{2-r+N}$, then 
\begin{multline}
\underset{i=1,\dots,r}{\sum_{0\le k_i\le N}}
\prod_{1\le i<j\le r}(q^{k_i}-q^{k_j})
  (1-aq^{k_i+k_j})\\
\times \prod_{i=1}^r\frac{(1-aq^{2k_i})(a,b,c,dx_i,e/x_i,q^{-N})_{k_i}}
{(1-a)(q,aq/b,aq/c,aq/dx_i,aqx_i/e,aq^{1+N})_{k_i}}\,
q^{k_i}\\
=\left(\frac {e^2}{ad}\right)^{\binom r2}
\prod_{1\le i<j\le r}(x_i-x_j)(1-dx_ix_j/e)
\prod_{i=1}^{r-1}\frac{(b,c,q^{-N})_i}{(aq^{2-r}/bc)_i}
\\
\times\prod_{i=1}^r
\frac{(aq,aq^{2-r}/bc)_N(aq^{2-r}/bdx_i,aq^{2-r}/cdx_i)_{N+1-r}}
{x_i^{2(r-1)}\,(aq^{2-i}/b,aq^{2-i}/c)_{N+1-i}(aq/dx_i,aq^{2-r}/bcdx_i)_{N}}.
\end{multline}
\end{Corollary}

The above $C_r$ Jackson summations all contain the factor \eqref{facdiff}
in the summand.
We now turn our attention to $C_r$ Jackson summations
of a different type, namely, of the type encountered in
\cite[Th.~4.2]{schlhypdet}. These also follow naturally from our $C_r$
nonterminating ${}_{10}\phi_9$ transformation in Corollary~\ref{cnnt109}.
For this case we put $gh=aq$, and $e_i=q^{-n_i}/x_i$, for $i=1,\dots,r$, in
\eqref{cnnt109gl}. On the right-hand side, because of $bf/\lambda=q^{1-r}$,
only the term with $S=\emptyset$ is non-zero. This term, similar
as in the proof of Corollary~\ref{cnt87}, is a determinant.
The result is the following (where we have replaced $f$ by $e$).

\begin{Corollary}[A $C_r$ Jackson summation]\label{cnt87c}
If $bc_id_ie=a^2q^{2-r+n_i}$, then 
\begin{multline}
\underset{i=1,\dots,r}{\sum_{0\le k_i\le n_i}}
\prod_{1\le i<j\le r}
\frac{(x_iq^{k_i}-x_jq^{k_j})(1-ax_ix_jq^{k_i+k_j})}
{(x_i-x_j)(1-ax_ix_j)}
\prod_{i=1}^r\frac{(1-ax_i^2q^{2k_i})}{(1-ax_i^2)}\\\times
\prod_{i=1}^r\frac{(ax_i^2,bx_i,c_ix_i,d_ix_i,ex_i,q^{-n_i})_{k_i}}
{(q,ax_iq/b,ax_iq/c_i,ax_iq/d_i,ax_iq/e,ax_i^2q^{1+n_i})_{k_i}}\,q^{k_i}\\
= (-1)^{\binom r2}b^{-\binom r2}q^{-2\binom r 3}
\prod_{1\le i<j\le r}\frac{1}{(x_i-x_j)(1-ax_ix_j)}\\
\times\prod_{i=1}^r\frac{(aqx_i^2,aq^{2-r}/bc_i,aq^{2-r}/bd_i)_{n_i}
(aq/c_id_i)_{n_i+i-r}(bx_i)_{r-1}}{(aqx_i/c_i,aqx_i/d_i,
aq^{2-r}/bc_id_ix_i)_{n_i}(aqx_i/b)_{n_i+1-r}(eq^{2+r-2i}/b)_{i-1}}
\\
\times \det_{1\le i,j\le r}\left(
\frac{(ex_i,ec_i/a,ed_i/a,eq^{-n_i}/ax_i)_{r-j}}
{(q^{2-r}/bx_i,q^{-n_i}ed_i/a,q^{-n_i}ec_i/a,aq^{2-r+n_i}x_i/b)_{r-j}}\right).
\end{multline}
\end{Corollary}

Analogously to Corollary~\ref{cnt87} we have two special cases
where the determinant evaluates as a product of linear factors,
giving rise to two different explicit Jackson summations.
The first one is $c_i=c$, $d_i=d$, $n_i=N$, $i=1,\dots,r$,
which gives the Tejasi sum in \cite[Th.~4.2]{schlhypdet}.
The case $x_i=x$, $d_i=d$,
$c_i=cq^{n_i}$,  $i=1,\dots,r$, gives back Corollary~\ref{cnt87a}.

\subsection{$\bmcr$ nonterminating $\bsixphifive$ summations}
We work out the various new $C_r$ extensions of Rogers' nonterminating
${}_6\phi_5$ summation (cf.\ \cite[Eq.~(2.7.1)]{grhyp})
following from
our results in Section~\ref{subsect87}.

First, we give the special case of Corollary~\ref{cnt87}
arising from formally replacing $b$ by $a^2q^{2-r+n_i}/c_id_ie_i$
(note that this is independent of $i$), and then letting
$n_i\to\infty$, for $i=1,\dots,r$. After subsequently
relabeling $e_i\mapsto b_i$, for $i=1,\dots,r$, we have the
following result.

\begin{Corollary}[A $C_r$ nonterminating ${}_6\phi_5$ summation]\label{cn65}
We have
\begin{multline}\label{cn65gl}
\sum_{k_1,\dots,k_r=0}^{\infty}
\prod_{1\le i<j\le r}(q^{k_i}-q^{k_j})
  (1-aq^{k_i+k_j})\\
\times \prod_{i=1}^r\frac{(1-aq^{2k_i})(a,b_i,c_i,d_i)_{k_i}}
{(1-a)(q,aq/b_i,aq/c_i,aq/d_i)_{k_i}}
\left(\frac{a^2q^{2-r}}{b_ic_id_i}\right)^{k_i}\\
=
q^{-\binom r3}
\det_{1\le i,j\le r}\left(\frac{(b_i,c_i,d_i)_{r-j}}
{(b_ic_id_i/a)_{r-j}}\right)
\prod_{i=1}^r\frac{(aq,aq/b_ic_i,aq/b_id_i,aq/c_id_i)_{\infty}}
{(aq/b_i,aq/c_i,aq/d_i,aq/b_ic_id_i)_{\infty}},
\end{multline}
provided $|aq^{2-r}/b_ic_id_i|<1$, for $i=1,\dots,r$.
\end{Corollary}

The following two results give cases where the determinant
on the right-hand side of \eqref{cn65gl} can be factored
by virtue of Lemma~\ref{lemdet1}.

For the first case, we specialize Corollary~\ref{cn65} by
putting $b_i=b$ and $c_i=c$, for $i=1,\dots,r$. 

\begin{Corollary}[A $C_r$ nonterminating ${}_6\phi_5$ summation]\label{cn65a}
We have
\begin{multline}
\sum_{k_1,\dots,k_r=0}^{\infty}
\prod_{1\le i<j\le r}(q^{k_i}-q^{k_j})
  (1-aq^{k_i+k_j})\\
\times \prod_{i=1}^r\frac{(1-aq^{2k_i})(a,b,c,d_i)_{k_i}}
{(1-a)(q,aq/b,aq/c,aq/d_i)_{k_i}}
\left(\frac{a^2q^{2-r}}{bcd_i}\right)^{k_i}\\
=
\prod_{1\le i<j\le r}(d_j-d_i)
\prod_{i=1}^r\frac{(b,c,bc/a)_{i-1}}{(bcd_i/a)_{r-1}}
\frac{(aq,aq/bc,aq/bd_i,aq/cd_i)_{\infty}}
{(aq/b,aq/c,aq/d_i,aq/bcd_i)_{\infty}},
\end{multline}
provided $|aq^{2-r}/bcd_i|<1$, for $i=1,\dots,r$.
\end{Corollary}

The other case results from choosing
$b_i=b$, $c_i=cx_i$, and $d_i=d/x_i$, for $i=1,\dots,r$,
in Corollary~\ref{cn65}.

\begin{Corollary}[A $C_r$ nonterminating ${}_6\phi_5$ summation]\label{cn65b}
We have
\begin{multline}
\sum_{k_1,\dots,k_r=0}^{\infty}
\prod_{1\le i<j\le r}(q^{k_i}-q^{k_j})
  (1-aq^{k_i+k_j})\\
\times \prod_{i=1}^r\frac{(1-aq^{2k_i})(a,b,cx_i,d/x_i)_{k_i}}
{(1-a)(q,aq/b,aq/cx_i,ax_iq/d)_{k_i}}
\left(\frac{a^2q^{2-r}}{bcd}\right)^{k_i}\\
=
c^{\binom r2}\prod_{1\le i<j\le r}(x_j-x_i)(1-d/cx_ix_j)\\\times
\prod_{i=1}^r\frac{(b)_{i-1}}{(bcd/a)_{i-1}}
\frac{(aq,aq/bcx_i,ax_iq/bd,aq/cd)_{\infty}}
{(aq/b,aq/cx_i,ax_iq/d,aq/bcd)_{\infty}},
\end{multline}
provided $|aq^{2-r}/bcd|<1$, for $i=1,\dots,r$.
\end{Corollary}

\subsection{$\bmcr$ terminating $\bsixphifive$ summations}
We now list the terminating versions of the above multivariable
${}_6\phi_5$ summations. These reduce to \cite[Eq.~(2.4.2)]{grhyp}
for $r=1$.

We start with the $d_i=q^{-n_i}$, $i=1,\dots,r$, case of Corollary~\ref{cn65}.

\begin{Corollary}[A $C_r$ terminating ${}_6\phi_5$ summation]\label{cnt65}
We have
\begin{multline}\label{cnt65gl}
\underset{i=1,\dots,r}{\sum_{0\le k_i\le n_i}}
\prod_{1\le i<j\le r}(q^{k_i}-q^{k_j})
  (1-aq^{k_i+k_j})\\
\times \prod_{i=1}^r\frac{(1-aq^{2k_i})(a,b_i,c_i,q^{-n_i})_{k_i}}
{(1-a)(q,aq/b_i,aq/c_i,aq^{1+n_i})_{k_i}}
\left(\frac{a^2q^{2-r+n_i}}{b_ic_i}\right)^{k_i}\\
=
q^{-\binom r3}
\det_{1\le i,j\le r}\left(\frac{(b_i,c_i,q^{-n_i})_{r-j}}
{(b_ic_iq^{-n_i}/a)_{r-j}}\right)
\prod_{i=1}^r\frac{(aq,aq/b_ic_i)_{n_i}}
{(aq/b_i,aq/c_i)_{n_i}}.
\end{multline}
\end{Corollary}

We give four different special cases of  Corollary~\ref{cnt65}, in each case
for which the determinant factors by virtue of Lemma~\ref{lemdet1}.

For the first case we put $b_i=b$ and $c_i=c$, for $i=1,\dots,r$.

\begin{Corollary}[A $C_r$ terminating ${}_6\phi_5$ summation]\label{cnt65a}
We have
\begin{multline}
\underset{i=1,\dots,r}{\sum_{0\le k_i\le n_i}}
\prod_{1\le i<j\le r}(q^{k_i}-q^{k_j})
  (1-aq^{k_i+k_j})\\
\times \prod_{i=1}^r\frac{(1-aq^{2k_i})(a,b,c,q^{-n_i})_{k_i}}
{(1-a)(q,aq/b,aq/c,aq^{1+n_i})_{k_i}}
\left(\frac{a^2q^{2-r+n_i}}{bc}\right)^{k_i}\\
=
\prod_{1\le i<j\le r}(q^{n_i}-q^{n_j})
\prod_{i=1}^r\frac{(b,c,bc/a)_{i-1}}
{(bc/a)_{r-1}}
\frac{(aq,aq^{2-r}/bc)_{n_i}}
{(aq/b,aq/c)_{n_i}}.
\end{multline}
\end{Corollary}

For the second case we put $b_i=b$ and $d_i=q^{-N}$, for $i=1,\dots,r$.

\begin{Corollary}[A $C_r$ terminating ${}_6\phi_5$ summation]\label{cnt65b}
We have
\begin{multline}
\underset{i=1,\dots,r}{\sum_{0\le k_i\le N}}
\prod_{1\le i<j\le r}(q^{k_i}-q^{k_j})
  (1-aq^{k_i+k_j})\\
\times \prod_{i=1}^r\frac{(1-aq^{2k_i})(a,b,c_i,q^{-N})_{k_i}}
{(1-a)(q,aq/b,aq/c_i,aq^{1+N})_{k_i}}
\left(\frac{a^2q^{2-r+N}}{bc_i}\right)^{k_i}\\
=
\prod_{1\le i<j\le r}(c_j-c_i)
\prod_{i=1}^r\frac{(b,q^{-N},bq^{-N}/a)_{i-1}}
{(bc_iq^{-N}/a)_{r-1}}\frac{(aq,aq/bc_i)_{N}}
{(aq/b,aq/c_i)_{N}}.
\end{multline}
\end{Corollary}

Note that Corollaries~\ref{cnt65a} and \ref{cnt65b} are equivalent
since they can be derived from each other by applying a standard
polynomial argument; cf.\ \cite[Th.~4.2]{milnetf}.
The same situation appears with
Corollaries~\ref{cnt65c} and \ref{cnt65d} below.

We now specialize Corollary~\ref{cnt65} by putting $b_i=bx_i$,
$c_i=c/x_i$, and $d_i=q^{-N}$, for $i=1,\dots,r$.

\begin{Corollary}[A $C_r$ terminating ${}_6\phi_5$ summation]\label{cnt65c}
We have
\begin{multline}
\underset{i=1,\dots,r}{\sum_{0\le k_i\le N}}
\prod_{1\le i<j\le r}(q^{k_i}-q^{k_j})
  (1-aq^{k_i+k_j})\\
\times \prod_{i=1}^r\frac{(1-aq^{2k_i})(a,bx_i,c/x_i,q^{-N})_{k_i}}
{(1-a)(q,aq/bx_i,ax_iq/c,aq^{1+N})_{k_i}}
\left(\frac{a^2q^{2-r+N}}{bc}\right)^{k_i}\\
=b^{\binom r2}
\prod_{1\le i<j\le r}(x_j-x_i)(1-c/bx_ix_j)
\prod_{i=1}^r\frac{(q^{-N})_{i-1}}{(bcq^{-N}/a)_{i-1}}\frac{(aq,aq/bc)_{N}}
{(aq/bx_i,ax_iq/c)_{N}}.
\end{multline}
\end{Corollary}

The last case, which extends the $c\mapsto cq^n$ version of
\cite[Eq.~(2.4.2)]{grhyp} comes from putting $b_i=b$, $c_i=cq^{n_i}$,
and $d_i=q^{-n_i}$, for $i=1,\dots,r$, in Corollary~\ref{cnt65}.

\begin{Corollary}[A $C_r$ terminating ${}_6\phi_5$ summation]\label{cnt65d}
We have
\begin{multline}
\underset{i=1,\dots,r}{\sum_{0\le k_i\le n_i}}
\prod_{1\le i<j\le r}(q^{k_i}-q^{k_j})
  (1-aq^{k_i+k_j})\\
\times \prod_{i=1}^r\frac{(1-aq^{2k_i})(a,b,cq^{n_i},q^{-n_i})_{k_i}}
{(1-a)(q,aq/b,aq^{1-n_i}/c,aq^{1+n_i})_{k_i}}
\left(\frac{a^2q^{2-r}}{bc}\right)^{k_i}\\
=
\prod_{1\le i<j\le r}(q^{-n_j}-q^{-n_i})(1-cq^{n_i+n_j})
\prod_{i=1}^r\frac{(b)_{i-1}}{(bc/a)_{i-1}}\frac{(aq,aq^{1-n_i}/bc)_{n_i}}
{(aq/b,aq^{1-n_i}/c)_{n_i}}.
\end{multline}
\end{Corollary}

\subsection{$\bmar$ terminating $\bthreephitwo$ summations}
In the following, we derive some multivariable extensions
of the terminating balanced ${}_3\phi_2$ (or $q$-Pfaff--Saalsch\"utz)
summation (cf.\ \cite[Eq.~(1.7.2)]{grhyp}).

In the $C_r$ Jackson summation in Corollary~\ref{cnt87},
we first remove the dependency of the
parameters by replacing $e_i$ by $a^2q^{2-r+n_i}/bc_id_i$, for
$i=1,\dots,r$. Then we replace $b$ by $aq^{2-r}/b$ and let $a\to 0$.
We perform the substitutions $c_i\mapsto a_i$, $d_i\mapsto b_i$,
for $i=1,\dots,r$, and $b\mapsto c$ and obtain the following result.

\begin{Corollary}[An $A_{r-1}$ terminating ${}_3\phi_2$ summation]\label{an32}
We have
\begin{multline}\label{an32gl}
\underset{i=1,\dots,r}{\sum_{0\le k_i\le n_i}}
\prod_{1\le i<j\le r}(q^{k_i}-q^{k_j})
\prod_{i=1}^r\frac{(a_i,b_i,q^{-n_i})_{k_i}}
{(q,cq^{r-1},a_ib_iq^{1-n_i}/c)_{k_i}}\,q^{k_i}\\
=
c^{\binom r2}\prod_{i=1}^r\frac{(c)_{r-1}}{(cq^{2r-2i})_{i-1}}
\frac{(c/a_i,c/b_i)_{n_i}}
{(c,c/a_ib_i)_{n_i}}\\\times
\det_{1\le i,j\le r}\left(\big(a_ib_iq^{-n_i}\big)^{j-r}
\frac{(a_i,b_i,q^{-n_i})_{r-j}}
{(c/a_i,c/b_i,cq^{n_i})_{r-j}}\right).
\end{multline}
\end{Corollary}

Next, we give three special cases where the determinant appearing
on the right-hand side of \eqref{an32gl} factors.

For the first case we put $a_i=a$, $b_i=b$, for $i=1,\dots,r$.
We would now require a limit case of the determinant evaluation
in Lemma~\ref{lemdet1}, the one which we explicitly stated in
\eqref{lemdet1a}. Equivalently, we start with
Corollary~\ref{cnt87a}, replace $e$ by $a^2q^{2-r}/bcd$, and then
$d$ by $aq^{2-r}/d$. We take $a\to 0$ and simultaneously
substitute the variables $c\mapsto a$ and $d\mapsto c$.
The result is the following.

\begin{Corollary}[An $A_{r-1}$ terminating ${}_3\phi_2$ summation]\label{an32a}
We have
\begin{multline}
\underset{i=1,\dots,r}{\sum_{0\le k_i\le n_i}}
\prod_{1\le i<j\le r}(q^{k_i}-q^{k_j})
\prod_{i=1}^r\frac{(a,b,q^{-n_i})_{k_i}}
{(q,cq^{r-1},abq^{1-n_i}/c)_{k_i}}\,q^{k_i}\\
=
\left(\frac c{ab}\right)^{\binom r2}
\prod_{1\le i<j\le r}(q^{n_i}-q^{n_j})
\prod_{i=1}^r\frac{(a,b)_{i-1}}{(c/a,c/b)_{i-1}}
\frac{(c/a,c/b)_{n_i}}
{(cq^{r-1},c/ab)_{n_i}}.
\end{multline}
\end{Corollary}

For the second case we specialize Corollary~\ref{an32} by putting
$a_i=a$, $n_i=N$, for $i=1,\dots,r$.
We again would require \eqref{lemdet1a}. Equivalently, we start with
Corollary~\ref{cnt87b}, replace $e$ by $a^2q^{2-r+N}/bcd$, and then
$b$ by $aq^{2-r}/b$. We take $a\to 0$, and then simultaneously
substitute the variables $c\mapsto a$, $x_i\mapsto b_i/d$,
for $i=1,\dots,r$, and $b\mapsto c$. This yields
the following result.

\begin{Corollary}[An $A_{r-1}$ terminating ${}_3\phi_2$ summation]\label{an32b}
We have
\begin{multline}
\underset{i=1,\dots,r}{\sum_{0\le k_i\le N}}
\prod_{1\le i<j\le r}(q^{k_i}-q^{k_j})
\prod_{i=1}^r\frac{(a,b_i,q^{-N})_{k_i}}
{(q,cq^{r-1},ab_iq^{1-N}/c)_{k_i}}\,q^{k_i}\\
=
\Big(\frac caq^N\Big)^{\binom r2}
\prod_{1\le i<j\le r}(1/b_i-1/b_j)
\prod_{i=1}^r\frac{(c)_{r-1}(a,q^{-N})_{i-1}}{(c/b_i)_{r-1}(c/a,cq^N)_{i-1}}
\frac{(c/a,c/b_i)_{N}}
{(c,c/ab_i)_{N}}.
\end{multline}
\end{Corollary}

We could have also applied a standard polynomial
argument to Corollary~\ref{an32a}, in order to derive
Corollary~\ref{an32b}.

For the next case, we set $a_i=a$, $b_i=bq^{n_i}$, for $i=1,\dots,r$,
in Corollary~\ref{an32}. Equivalently, we start with
Corollary~\ref{cnt87a}, replace $d$ by $a^2q^{2-r}/bce$, and then
$b$ by $aq^{2-r}/b$. We take $a\to 0$, and then simultaneously
substitute the variables $b\mapsto c$, $c\mapsto a$ and $e\mapsto b$,
and arrive at the following result which reduces to the
$b\mapsto bq^n$ case of \cite[Eq.~(1.7.2)]{grhyp}.

\begin{Corollary}[An $A_{r-1}$ terminating ${}_3\phi_2$ summation]\label{an32c}
We have
\begin{multline}
\underset{i=1,\dots,r}{\sum_{0\le k_i\le n_i}}
\prod_{1\le i<j\le r}(q^{k_i}-q^{k_j})
\prod_{i=1}^r\frac{(a,bq^{n_i},q^{-n_i})_{k_i}}
{(q,cq^{r-1},abq/c)_{k_i}}\,q^{k_i}\\
=
q^{-\binom r3}a^{-\binom r2}
\prod_{1\le i<j\le r}(q^{-n_i}-q^{-n_j})(1-bq^{n_i+n_j})\\\times
\prod_{i=1}^r\frac{(a)_{i-1}}{(c/a,bq^{2-r}/c)_{i-1}}
\frac{(c/a,cq^{r-1-n_i}/b)_{n_i}}
{(cq^{r-1},cq^{-n_i}/ab)_{n_i}}.
\end{multline}
\end{Corollary}

For the third case, we choose $a_i=ax_i$, $b_i=b/x_i$, and 
$n_i=q^{-N}$ for $i=1,\dots,r$. Equivalently, we start with
Corollary~\ref{cnt87b}, replace $c$ by $a^2q^{2-r+N}/bde$, and then
$b$ by $aq^{2-r}/b$. We take $a\to 0$, and then simultaneously
substitute the variables $d\mapsto a$, $e\mapsto b$ and $b\mapsto c$,
and arrive at the following result.

\begin{Corollary}[An $A_{r-1}$ terminating ${}_3\phi_2$ summation]\label{an32d}
We have
\begin{multline}
\underset{i=1,\dots,r}{\sum_{0\le k_i\le N}}
\prod_{1\le i<j\le r}(q^{k_i}-q^{k_j})
\prod_{i=1}^r\frac{(ax_i,b/x_i,q^{-N})_{k_i}}
{(q,cq^{r-1},abq^{1-N}/c)_{k_i}}\,q^{k_i}\\
=
q^{\binom r3}\Big(\frac cbq^N\Big)^{\binom r2}
\prod_{1\le i<j\le r}(x_j-x_i)(1-b/ax_ix_j)\\\times
\prod_{i=1}^r\frac{(c)_{r-1}(c/ab,q^{-N})_{i-1}}
{(cq^N)_{i-1}(c/ax_i,cx_i/b)_{r-1}}
\frac{(c/ax_i,cx_i/b)_{N}}
{(c,c/ab)_{N}}.
\end{multline}
\end{Corollary}

Note that Corollaries~\ref{an32c} and \ref{an32d} are equivalent
since they can be derived from each other by applying a standard
polynomial argument.

\subsection{$\bmar$ nonterminating $\btwophione$ and $\bonephizero$
summations}
We briefly give some $A_{r-1}$ $q$-Gau{\ss} summations and
$A_{r-1}$ $q$-binomial theorems which follow from our results.
We omit writing out the terminating versions but we have made an
effort in presenting all the nonterminating ones explicitly.

We start with our $A_{r-1}$ $q$-Gau{\ss} summations.
If, in Corollary~\ref{an32},
we let $n_i\to\infty$, for all $i=1,\dots,r$, we obtain

\begin{Corollary}[An $A_{r-1}$ $q$-Gau{\ss} summation]\label{an21}
We have
\begin{multline}\label{an21gl}
\sum_{k_1,\dots,k_r=0}^{\infty}
\prod_{1\le i<j\le r}(q^{k_i}-q^{k_j})
\prod_{i=1}^r\frac{(a_i,b_i)_{k_i}}
{(q,cq^{r-1})_{k_i}}\left(\frac c{a_ib_i}\right)^{k_i}\\
=
(-c)^{\binom r2}q^{\binom r3}
\prod_{i=1}^r\frac{(c)_{r-1}}{(cq^{2r-2i})_{i-1}}
\frac{(c/a_i,c/b_i)_{\infty}}
{(c,c/a_ib_i)_{\infty}}
\det_{1\le i,j\le r}\left(\big(a_ib_i\big)^{j-r}
\frac{(a_i,b_i)_{r-j}}
{(c/a_i,c/b_i)_{r-j}}\right),
\end{multline}
where $|c/a_ib_i|<1$, for $i=1,\dots,r$.
\end{Corollary}

To evaluate the determinant in factored form we choose different
specializations of the parameters.

We can choose $a_i=ax_i$ and $b_i=b/x_i$, for $i=1,\dots,r$,
and then substitute $c\mapsto cq^{1-r}$,
which gives the following.

\begin{Corollary}[An $A_{r-1}$ $q$-Gau{\ss} summation]\label{an21a}
We have
\begin{multline}
\sum_{k_1,\dots,k_r=0}^{\infty}
\prod_{1\le i<j\le r}(q^{k_i}-q^{k_j})
\prod_{i=1}^r\frac{(ax_i,b/x_i)_{k_i}}
{(q,c)_{k_i}}\left(\frac{cq^{1-r}}{ab}\right)^{k_i}\\
=
\Big(\frac c{bq}\Big)^{\binom r2}q^{-\binom r3}
\prod_{1\le i<j\le r}(x_i-x_j)(1-b/ax_ix_j)
\prod_{i=1}^r\frac{(c/ax_i,cx_i/b)_{\infty}}
{(c,cq^{1-i}/ab)_{\infty}},
\end{multline}
where $|cq^{1-r}|<1$.
\end{Corollary}

Here is the case $a_i=a$, for $i=1,\dots,r$, of Corollary~\ref{an21},
with $c\mapsto cq^{1-r}$.

\begin{Corollary}[An $A_{r-1}$ $q$-Gau{\ss} summation]\label{an21b}
We have
\begin{multline}
\sum_{k_1,\dots,k_r=0}^{\infty}
\prod_{1\le i<j\le r}(q^{k_i}-q^{k_j})
\prod_{i=1}^r\frac{(a,b_i)_{k_i}}
{(q,c)_{k_i}}\left(\frac{cq^{1-r}}{ab_i}\right)^{k_i}\\
=
\Big(\frac c{aq}\Big)^{\binom r2}q^{-2\binom r3}
\prod_{1\le i<j\le r}(1/b_i-1/b_j)
\prod_{i=1}^r(a)_{i-1}
\frac{(cq^{i-1}/a,c/b_i)_{\infty}}
{(c,cq^{1-r}/ab_i)_{\infty}},
\end{multline}
where $|cq^{1-r}/ab_i|<1$, for $i=1,\dots,r$.
\end{Corollary}

Next, we give some $A_{r-1}$ $q$-binomial theorems.
If in Corollary~\ref{an21} we replace $b_i$ by $c/a_iz_i$, for $i=1,\dots,r$,
and then let $c\to 0$, we obtain the following summation.

\begin{Corollary}[An $A_{r-1}$ $q$-binomial theorem]\label{an10}
We have
\begin{multline}\label{an10gl}
\sum_{k_1,\dots,k_r=0}^{\infty}
\prod_{1\le i<j\le r}(q^{k_i}-q^{k_j})
\prod_{i=1}^r\frac{(a_i)_{k_i}}
{(q)_{k_i}}z_i^{k_i}\\
=
(-1)^{\binom r2}q^{\binom r3}
\det_{1\le i,j\le r}\left(z_i^{r-j}
\frac{(a_i)_{r-j}}
{(a_iz_i)_{r-j}}\right)
\prod_{i=1}^r\frac{(a_iz_i)_{\infty}}
{(z_i)_{\infty}},
\end{multline}
where $|z_i|<1$, for $i=1,\dots,r$.
\end{Corollary}

Note that if we let $b_i=q$ for $i=1,\dots,r$ in
Proposition~\ref{an11}, we alternatively get the identity
\begin{equation}\label{an10vgl}
\sum_{k_1,\dots,k_r=0}^{\infty}
\prod_{1\le i<j\le r}(q^{k_i}-q^{k_j})
\prod_{i=1}^r\frac{(a_i)_{k_i}}
{(q)_{k_i}}z_i^{k_i}
=
\det_{1\le i,j\le r}\left(
\frac{(z_i)_{r-j}}
{(a_iz_i)_{r-j}}\right)
\prod_{i=1}^r\frac{(a_iz_i)_{\infty}}
{(z_i)_{\infty}}.
\end{equation}
In fact, comparing the two right-hand sides of \eqref{an10gl} and
\eqref{an10vgl}, we observe that the following transformation of
determinants must hold:
\begin{equation}\label{dettf}
\det_{1\le i,j\le r}\left(
\frac{(z_i)_{r-j}}
{(a_iz_i)_{r-j}}\right)
=
(-1)^{\binom r2}q^{\binom r3}
\det_{1\le i,j\le r}\left(z_i^{r-j}
\frac{(a_i)_{r-j}}
{(a_iz_i)_{r-j}}\right).
\end{equation}
This determinant transformation can also be easily proved in a
more natural way. (This has been kindly communicated to us by
Christian Krattenthaler.) Note that by reversing the order of
columns of the matrix on the left hand side of \eqref{dettf}
(by which the determinant gets multiplied with $(-1)^{\binom r2}$),
we have the matrix $\big((z_i)_{j-1}/(a_iz_i)_{j-1}\big)_{1\le i,j,\le r}$.
If this matrix is now being multiplied with the lower-triangular
matrix $\big(q^{(j-1)(1-k)}(q^{1-k})_{k-j}/(q)_{k-j}\big)_{1\le j,k\le r}$
(which has determinant $q^{\binom r2-2\binom r3}$),
we obtain, after application of the $q$-Chu--Vandermonde summation
\cite[Eq.~(II.6)]{grhyp} the matrix $\big((-1)^{k-1}z_i^{k-1}q^{-\binom k2}
(a_i)_{k-1}/(a_iz_i)_{k-1}\big)_{1\le i,k\le r}$.
Changing back the order of columns and taking determinants
we immediately establish \eqref{dettf}.

We complete our listing of summations by giving three special cases
of Corollary~\ref{an10}
for which the determinant
factors by virtue of Lemma~\ref{lemdet1}.

The first case is $a_i=a$, for $i=1,\dots,r$.

\begin{Corollary}[An $A_{r-1}$ $q$-binomial theorem]\label{cn10a}
We have
\begin{multline}
\sum_{k_1,\dots,k_r=0}^{\infty}
\prod_{1\le i<j\le r}(q^{k_i}-q^{k_j})
\prod_{i=1}^r\frac{(a)_{k_i}}
{(q)_{k_i}}z_i^{k_i}\\
=
q^{\binom r3}\prod_{1\le i<j\le r}(z_j-z_i)
\prod_{i=1}^r(a)_{i-1}
\frac{(az_iq^{r-1})_{\infty}}{(z_i)_{\infty}},
\end{multline}
where $|z_i|<1$, for $i=1,\dots,r$.
\end{Corollary}

The second case is $a_i=a/x_i$, and $z_i=zx_i$ for $i=1,\dots,r$.

\begin{Corollary}[An $A_{r-1}$ $q$-binomial theorem]\label{cn10b}
We have
\begin{multline}
\sum_{k_1,\dots,k_r=0}^{\infty}
\prod_{1\le i<j\le r}(q^{k_i}-q^{k_j})
\prod_{i=1}^r\frac{(a/x_i)_{k_i}}
{(q)_{k_i}}(zx_i)^{k_i}\\
=
z^{\binom r2}q^{\binom r3}\prod_{1\le i<j\le r}(x_j-x_i)
\prod_{i=1}^r
\frac{(azq^{i-1})_{\infty}}{(zx_i)_{\infty}},
\end{multline}
where $|zx_i|<1$, for $i=1,\dots,r$.
\end{Corollary}

The third case is $z_i=z$ for $i=1,\dots,r$.

\begin{Corollary}[An $A_{r-1}$ $q$-binomial theorem]\label{cn10c}
We have
\begin{multline}
\sum_{k_1,\dots,k_r=0}^{\infty}
\prod_{1\le i<j\le r}(q^{k_i}-q^{k_j})
\prod_{i=1}^r\frac{(a_i)_{k_i}}
{(q)_{k_i}}z^{k_i}\\
=
z^{\binom r2}q^{2\binom r3}\prod_{1\le i<j\le r}(a_i-a_j)
\prod_{i=1}^r
\frac{(a_izq^{r-1})_{\infty}}{(zq^{i-1})_{\infty}},
\end{multline}
where $|z|<1$.
\end{Corollary}

\section{Elliptic extensions}\label{secell}

In this section we extend some of our results to the case of 
elliptic hypergeo\-metric series. Note that
these satisfy fewer identities than basic
hypergeometric series. Roughly speaking, the reason is that
elliptic hypergeometric identities generalize Riemann's
addition formula for theta functions, which is more complicated than the
addition formulas for trigonometric functions.
Moreover, infinite elliptic
hypergeometric series lead to serious problems of convergence; 
in particular, no elliptic
analogue of \eqref{109ntgl} is known. This means that  we can only  extend
the results of Sections \ref{sect109} and \ref{subsect87} to the elliptic case.
We refer to \cite{spths} for a detailed discussion of the balanced and
well-poised conditions for elliptic hypergeometric series, and their
relation to  modular invariance of the series.

Our elliptic extensions involve a fixed parameter
$p$ such that $|p|<1$. We write $\theta(x):=(x,p/x;p)_\infty$ and 
 define elliptic shifted factorials by
\begin{equation}\label{defesf}(a;q,p)_k:=\prod_{j=0}^{k-1}\theta(aq^j).
\end{equation}
When $p=0$, $\theta(x)=1-x$ and we recover the $q$-shifted factorials
used before. Since $p$ and $q$ are fixed we omit them from the
notation, writing
$$(a)_k:=(a;q,p)_k.$$
To use the same shorthand notation as in \eqref{defqsf}  might 
seem confusing, but has the advantage that we can almost 
use our previous results as they stand, just interpreting
the symbol $(a)_k$ differently. The only other change we have to make is
that all factors of the form  $x-y$ should be replaced by
$x\,\theta(y/x)=-y\,\theta(x/y)$. Thus, the ubiquitous factor
$$\prod_{1\leq i<j\leq r}(t_i-t_j)(1-at_it_j)\prod_{i=1}^r(1-at_i^2) $$
is replaced by
$$\prod_{1\leq i<j\leq r}t_i\,\theta(t_j/t_i)\theta(at_it_j)
\prod_{i=1}^r\theta(at_i^2).$$
All other factors, such as $q^{k_i}$ and and $(\lambda/a)^{\binom r2}$
in Corollary \ref{cnt109gl}, are left untouched.

\begin{Theorem}
\label{thell}
The following results have elliptic analogues, given as explained above:
\emph{Corollary \ref{cnt109gl}, Corollary \ref{cnt109n}, Corollary
\ref{cnt109i2}, Corollary \ref{cnt87}, Corollary \ref{cnt87a},
Corollary \ref{cnt87b}} and \emph{Corollary \ref{cnt87c}}.
\end{Theorem}

When $r=1$, Theorem \ref{thell} reduces to the elliptic Jackson
summation and the elliptic Bailey transformation of Frenkel and Turaev
\cite{ft}. 

\begin{proof}
First note that to prove Corollary \ref{cnt109gl} we did not use the
general case of \eqref{109ntgl}. We only used the identity obtained
after multiplying by $(bd/a)_\infty$ and then letting $bd/a=q^{-N}$.
This is the terminating Bailey transformation
\eqref{109gl}, whose elliptic analogue is known to hold \cite{ft}. 
We also needed the determinant evaluation of Lemma \ref{lemdet1}, whose
elliptic analogue was obtained by Warnaar \cite{warnell}. Apart from
these fundamental results, we only used elementary identities for
$q$-shifted factorials that also hold in the elliptic case (indeed,
they would hold if $\theta(x)$ in \eqref{defesf} was replaced by any function
satisfying $\theta(1/x)=-\theta(x)/x$). Thus, the proof of Corollary 
 \ref{cnt109gl} extends immediately to the elliptic case.

The derivations of Corollary \ref{cnt109n}, Corollary \ref{cnt87},
Corollary  \ref{cnt87a} and Corollary \ref{cnt87b} from Corollary 
 \ref{cnt109gl} involve only elementary identities and Lemma \ref{lemdet1},
and thus extend to the elliptic case.

In the proof of Corollary \ref{cnt109i2} we made essential use of
infinite series, so we need an alternative approach. 
 We  sketch a way  to derive the elliptic case of
Corollary  \ref{cnt109i2} from the elliptic case of
Corollary \ref{cnt109n}. We use an elliptic version of the
``polynomial argument''; this idea also occurs in 
 \cite{r2} and  \cite{warnell}.
We start with the elliptic Corollary \ref{cnt109n},
multiply it by $\prod_{i=1}^r(aq^{1-n_i}/e)_{n_i}$ and then let
$e\rightarrow aq^{-N}$ with $N$ a non-negative integer. 
On the left-hand side we then have the factor
$$\frac{(aq^{1-n_i}/e)_{n_i}}{(aq^{1-n_i}/e)_{k_i}}
\rightarrow(q^{1+k_i+N-n_i})_{n_i-k_i},
$$
which vanishes for $n_i-k_i\geq N+1$. On the right-hand side, we have
$$(eq^{1-r}/a)_{k_i}=(q^{1-r-N})_{k_i} $$
which vanishes for $k_i\geq N+r$. Thus,
after  replacing $k_i$ by
$n_i-k_i$ on the left-hand side, we obtain a transformation of the form
$$\underset{i=1,\dots,r}{\sum_{0\le k_i\le \min(n_i,N)}}(\dotsm)=
\underset{i=1,\dots,r}{\sum_{0\le k_i\le \min(n_i,N+r-1)}}(\dotsm).$$
A straightforward computation  reveals that this  is
equivalent to the case $ex_i=q^{-n_i}$ of the elliptic analogue of
Corollary  \ref{cnt109i2}.

Now let $f(x_1,\dots,x_n)$ denote the  left-hand minus the right-hand side
of the elliptic analogue of \eqref{cnt109i2gl}. 
We have proved that $f$ vanishes if $x_i=q^{-n_i}/e$ for all $i$. A 
computation, using $\theta(px)=-x^{-1}\theta (x)$ and
$(pa)_k=(-1)^kq^{-\binom k2}a^{-k}(a)_k$,
shows that 
$$f(x_1,\dots,px_i,\dots,x_n)= f(x_1,\dots,x_n).$$
Thus, $f$ vanishes also when $x_i=p^kq^{-n_i}/e$ with $k$ an integer.
For generic values of the parameters, these zeroes have a limit point
in which $f$ is analytic, so $f$  is identically zero by analytic
continuation. Finally, 
analytic continuation in the remaining parameters extends the
identity to non-generic situations.

It remains to treat
Corollary \ref{cnt87c}. In this case, we only used \eqref{109ntgl} in
the case when $e=q^{-l}$ and $bf/\lambda=q^{-m}$ for non-negative
integers $l$, $m$. This is a
transformation between two finite sums. If in addition
$b=q^{-n}$, it reduces to a special case of \eqref{109gl}. In the
elliptic case, one may then use the same argument as above to prove the
corresponding identity. That is, one first proves that the (known) case 
$b=q^{-n}$ implies the case $b=p^kq^{-n}$ for $k\in\mathbb Z$,
and then extends it to general
$b$ by analytic continuation. Once the needed transformation formula 
is known,  the proof of Corollary \ref{cnt87c} extends immediately
to the elliptic case.
\end{proof}

\end{document}